\documentclass[11pt,a4paper,reqno]{amsart}
\usepackage[applemac]{inputenc}
\usepackage[T1]{fontenc}
\usepackage{amsmath}
\usepackage{amsthm}
\usepackage{amsfonts}
\usepackage{amssymb}
\usepackage{graphicx}
\usepackage{amsbsy}
\usepackage{mathrsfs}
\usepackage{bbm}
\newtheorem{theorem}{Theorem}[section]
\newtheorem{lemma}[theorem]{Lemma}

\newtheorem{proposition}[theorem]{Proposition}

\theoremstyle{remark}
\newtheorem{remark}[theorem]{\it \bf{Remark}\/}

\numberwithin{equation}{section}
\catcode`@=11
\def\section{\@startsection{section}{1}%
  \z@{1.5\linespacing\@plus\linespacing}{.5\linespacing}%
  {\normalfont\bfseries\large\centering}}
\catcode`@=12
\newcommand{\be}{\begin{equation}}
\newcommand{\ee}{\end{equation}}
\newcommand{\bea}{\begin{eqnarray}}
\newcommand{\eea}{\end{eqnarray}}
\newcommand{\bee}{\begin{eqnarray*}}
\newcommand{\eee}{\end{eqnarray*}}

\def\pa{\partial}

\def\RR{\mathbb{R}}

\def\fref#1{{\rm (\ref{#1})}}

\def\G{{\Gamma}}

\catcode`@=11
\def\supess{\mathop{\operator@font Sup\,ess}}
\catcode`@=12

\def\bt{\tilde{b}}

\def\RR{\mathbb{R}}

\def\e{\varepsilon}

\def\g#1{{\bf #1}}

\def\fref#1{{\rm (\ref{#1})}}

\def\R2+{\RR ^2_+}

\def\lsl{\frac{\lambda_s}{\lambda}}

\def\pa{\partial}

\def\lim{\mathop{\rm lim}}

\def\l{\lambda}

\def\log{{\rm log}}

\def\et{\tilde{\e}}

\def\lsl{\frac{\lambda_s}{\lambda}}

\def\Psih{\hat{\Psi}}

\def\ut{\tilde{u}}

\def\qbt{\tilde{Q}_b}

\def\tt{\tilde{T}}

\def\pa{\partial}

\def\et{\tilde{\e}}

\def\eh{\hat{\e}}

\def\alphah{\hat{\alpha}}

\def\pa{\partial}
\def\alphat{\tilde{\alpha}}
\def\Psit{\tilde{\Psi}}
\def\qbh{\hat{Q}_b}
\def\wh{\hat{w}}
\def\zh{\hat{z}}
\def\vt{\tilde{v}}

\title[]{Stable blow up dynamics for the 1-corotational energy critical harmonic heat flow}
\author[P. Rapha\"el]{Pierre Rapha\"el}
\address{Institut de Math\'ematiques de Toulouse, Universit\'e Paul  Sabatier, Toulouse, France}
\email{pierre.raphael@math.univ-toulouse.fr}
\author[R. Schweyer]{R\'emi Schweyer}
\address{Institut de Math\'ematiques de Toulouse, Universit\'e Paul  Sabatier, Toulouse, France}
\email{remi.schweyer@math.univ-toulouse.fr}

\begin{document}
\maketitle

\begin{abstract}
We exhibit a stable finite time blow up regime for the 1-corotational energy critical harmonic heat flow from $\Bbb R^2$ into a smooth compact revolution surface of $\Bbb R^3$ which reduces to the semilinear parabolic problem $$\partial_t u -\pa^2_{r} u-\frac{\pa_r u}{r} + \frac{f(u)}{r^2}=0$$ for a suitable class of functions $f$. The corresponding initial data can be chosen smooth, well localized and arbitrarily close to the ground state harmonic map in the energy critical topology. We give sharp asymptotics on the corresponding singularity formation which occurs through the concentration of a universal bubble of energy at the speed predicted in \cite{heatflow}. Our approach lies in the continuation of the study of the 1-equivariant energy critical wave map and Schr\"odinger map with $\Bbb S^2$ target in \cite{RR}, \cite{MRR}.
\end{abstract}


\section{Introduction}


\subsection{Setting of the problem}


The harmonic heat flow between two embedded Riemanian manifolds $(N,g_N),(M,g_M)$ is the gradient flow associated to the Dirichlet energy of maps from $N\to M$: 
\be
\label{hfgeneral}
\left\{\begin{array}{ll}\pa_tv={\Bbb P}_{T_vM}(\Delta_{g_N}v)\\
v_{|t=0}=v_0\end{array}\right . \ \ (t,x)\in \Bbb R\times N, \ \ v(t,x)\in M
\ee where ${\Bbb P}_{T_vM}$ is the projection onto the tangent space to $M$ at v. The special case $N=\Bbb R^2$, $M=\Bbb S^2$ corresponds to the harmonic heat flow to the 2-sphere
\be
\label{harmonicheatflow}
\pa_tv=\Delta v+|\nabla v|^2v, \ \ (t,x)\in \Bbb R\times \Bbb R^2,\ \ v(t,x)\in \Bbb S^2
\ee
which appears in cristal physics and is related to the Landau Lifschitz equation of ferromagnetism, we refer to \cite{heatflow}, \cite{matano}, \cite{NT1}, \cite{NT2} and references therein for a complete introduction to this class of problems. We shall from now on restrict our discussion to the case: $$N=\Bbb R^2.$$
Local existence of solutions emanating from smooth data is well known. Note that the Dirichlet energy is dissipated by the flow $$\frac{d}{dt}\left\{\int_{\Bbb R^2}|\nabla v|^2\right\}=-2\int_{\Bbb R^2}|\pa_tv|^2$$ and left invariant by the scaling symmetry $$u_\l(t,x)=u(\lambda^2t, \lambda x).$$ Hence the problem is {\it energy critical} and a singularity formation by energy concentration is possible. By the works of Struwe \cite{Struwe}, Ding and Tian \cite{dingtian}, Qing and Tian \cite{qingtian} (see Topping \cite{Topping} for a complete history of the problem), it is known that if occuring, concentration implies the bubbling off of a non trivial harmonic map at a finite number of blow up points 
\be
\label{cnknneooe}
v(t_i,a_i+\lambda(t_i)x)\to  Q_i, \ \ \l(t_i)\to 0
\ee
locally in space. In particular, this shows global existence on negatively curved target where no nontrivial harmonic map exists.


\subsection{Corotational flows}


The existence of blow up solutions has been proved in various different geometrical settings, see in particular Chang, Ding, Ye \cite{CDY}, Coron and Ghidaglia \cite{CD}, Qing and Tian \cite{qingtian}, Topping \cite{Topping}. We shall restrict in this paper onto flows with symmetries which are better understood.\\
Let a smooth closed curve in the plane parametrized by arclength $$u\in[-\pi,\pi]\mapsto \left|\begin{array}{ll}g(u)\\z(u)\end{array}\right., \ \ (g')^2+(z')^2=1,$$ where 
\be
\label{assumtiong}
(H)\ \ \left\{\begin{array}{lll} g\in \mathcal C^{\infty}(\Bbb R) \ \ \mbox{is odd and $2\pi$ periodic},\\
g(0)=g(\pi)=0, \ \ g(u)>0\ \ \mbox{for}\ \ 0<u<\pi,\\
g'(0)=1, \ \ g'(\pi)=-1,\end{array}\right.
\ee
then the revolution surface $M$ with parametrization $$(\theta,u)\in[0,2\pi]\times[0,\pi]\mapsto \left|\begin{array}{lll} g(u)\cos\theta\\g(u)\sin\theta\\z(u)\end{array}\right.,$$ is a smooth\footnote{see eg  \cite{GHL}} compact revolution surface of $\Bbb R^3$ with metric $(du)^2+(g(u))^2(d\theta)^2$. Given a homotopy degree $k\in \Bbb Z^*$, the k-corotational reduction to \fref{hfgeneral} corresponds to solutions of the form 
\be
\label{corotshpoere} 
v(t,r)= \left|\begin{array}{lll} g(u(t,r))\cos(k\theta)\\g(u(t,r))\sin(k\theta)\\z(u(t,r))\end{array}\right.
\ee which leads to the semilinear parabolic equation\footnote{see \fref{formuleprojection}}: 
\be
\label{mapg}
\left\{\begin{array}{ll}\pa_tu-\pa^2_ru-\frac{\pa_ru}{r}+k^2\frac{f(u)}{r^2}=0,\\
u_{t=0}=u_0\end{array}\right.\ \ f=gg'.
\ee
The k-corotational Dirichlet energy becomes 
\be
\label{diricheltenergy}
E(u)=\int_0^{+\infty}\left[|\pa_ru|^2+k^2\frac{(g(u))^2}{r^2}\right] rdr
\ee
and is minimized along maps with boundary conditions 
\be
\label{boundary}
u(0)=0, \ \ \lim_{r\to+\infty}u(r)=\pi
\ee
 onto the least harmonic map $Q_k$ which is the unique -up to scaling- solution to
 \be
 \label{eqgroudnsttae}
 r\pa_rQ_k=g(Q_k)
 \ee satisfying \fref{boundary}, see for example \cite{cote}.\\
In the case of $\Bbb S^2$ target $g(u)=\sin u$, the harmonic map is explicitely given by $$Q_k(r)=2\tan^{-1}(r^k).$$ In the series of works by Guan, Gustaffson,  Tsai \cite{NT1}, Gustaffson, Nakanishi, Tsai \cite{NT2}, $Q_k$ is proved to be {\it stable} by the flow \fref{mapg} for $k\geq 3$, and in particular no blow up will occur near $Q_k$. Moreover, eternally oscillating solutions and infinite time blow up solutions are exhibited for $k=2$. For the degree $k=1$ least energy harmonic map $Q\equiv Q_1$ with $\Bbb D^2$ initial manifold and $\Bbb S^2$ target, the formal analysis by Van den Bergh, Hulshof and King \cite{heatflow} suggests through matching asymptotics the existence of a {\it stable generic} blow up regime with $$u(t,r)\sim Q\left(\frac{r}{\lambda(t)}\right), \ \ \lambda(t)\sim \frac{T-t}{|\log (T-t)|^2}.$$ In this direction and for $k=1$, Angenent, Hulshof and Matano exhibit in \cite{matano} a class of corotational solutions which blow up in finite time with an estimate: $$\l(t)=o(T-t)\ \ \mbox{as} \ \ t\to T.$$The maximum principle plays an important role in this analysis. The sharp description of the singularity formation for $k=1$ and in particular the understanding of the generic regime thus remain open.\\ 
More generally, let us recall that the derivation of the blow up speed for energy critical parabolic problems is poorly understood, and for example the derivation of sharp asymptotics of type II blow up for the energy critical semilinear problem $$\pa_tu=\Delta u+u^{\frac{N+2}{N-2}}, \ \ (t,x)\in \Bbb R\times \Bbb R^N, \ \ N\geq 3$$ is open.


\subsection{Statement of the result}


Let $  \mathcal Q$ the least energy harmonic map with degree $1$ generated by the $Q\equiv Q_1$ solution to \fref{eqgroudnsttae}, explicitely:
\be
\label{defharmmap}
\mathcal Q(x)=\left|\begin{array}{lll} g(Q(t,r))\cos\theta\\g(Q(t,r))\sin\theta\\z(Q(t,r))\end{array}\right.
\ee
For an integer $i\geq 1$, we let $\dot{H}^i$ be the completion of $\mathcal C^{\infty}_c(\Bbb R^2,\Bbb R^3)$ for the norm $$\|v\|_{\dot{H}^i}=\|\Delta^{\frac{i}{2}} v\|_{L^2}.$$ 
We claim the existence and stability of a universal blow up regime emerging from 1-equivariant smooth data arbitrarily close to $\mathcal Q$ in the energy critical topology, together with sharp asymptotics on the singularity formation.

\begin{theorem}[Stable blow up dynamics for the 1-corotational heat flow]
\label{thmmain}
Let $k=1$ and $g$ satisfy \fref{assumtiong}. Let $\mathcal Q$ be the least energy harmonic map given by \fref{defharmmap}.  Then there exists an open set $\mathcal O$ of 1-corotational initial data of the form $$v_0=\mathcal Q+\e_0, \ \ \e_0\in \mathcal O\subset \dot{H}^1\cap\dot{H}^4$$ such that the corresponding solution $v\in \mathcal C([0,T),\dot{H}^1\cap\dot{H}^4)$ to \fref{hfgeneral} blows up in finite time $0<T=T(u_0)<+\infty$ according to the following universal scenario:\\
{\em (i) Universality of the concentrating bubble}: there exists an asympotic profile $v^*\in \dot{H}^1$ and $\lambda\in \mathcal C^1([0,T),\Bbb R^*_+)$ such that 
\be
\label{convustarb}
\lim_{t\to T}\left\|v(t,x)-\mathcal Q\left(\frac{x}{\lambda(t)}\right)-v^*\right\|_{\dot H^1}=0.
\ee
{\em (ii) Sharp asymptotics}: the blow up speed is given by
\be
\label{universallawkgeq}
\lambda(t)=c(v_0)(1+o(1))\frac{T-t}{|\log(T-t)|^2} \ \ \mbox{as} \ \ t\to T
\ee
for some $c(v_0)>0$.\\
{\em (iii) Regularity of the asymptotic profile}: there holds the additional regularity
\be
\label{regularityustar}
v^*\in \dot{H}^2.
\ee
\end{theorem}

In other words, there exists a generic blow up regime with the law \fref{universallawkgeq} as predicted in \cite{heatflow} for $g(u)=\sin u$, and blow up in this regime occurs by the concentration of a universal and quantized bubble of energy.\\

{\it Comments on the result:}\\

{\it 1. Energy method}: Following the strategy developped in \cite{MR1}, \cite{MR4}, \cite{RDuke}, \cite{RR}, \cite{MRR}, our strategy of proof proceeds first through the construction of suitable approximate solutions, and then the control of the remainding radiation through a robust {\it energy method}. In particular, we make no use of the maximum principle, and hence we expect our strategy to be applicable to more complicated parabolic {\it systems} among which the full problem \fref{hfgeneral}. Note also that parabolic problems cannot be solved backwards in time and involve smooth data. In this sense the construction of blow up solutions requires to follow the flow of smooth solutions forward in time and cannot be achieved by solving from blow up time for rough data as in \cite{Mm}, \cite{KST}. The set of initial data we construct in the proof of Theorem \ref{thmmain} contains compactly supported $C^{\infty}$ 1-corotational functions.\\

{\it 2. Regularity of the asymptotic profile}: The regularity of the asymptotic profile \fref{regularityustar} is a completely new feature with respect to the regularity obtained in \cite{MR5}, \cite{RR} where the profile is just in the critical space. This would also allow one to quantify the convergence rate \fref{convustarb} and bound the error polynomially in time, which is  a substantial improvement on the general convergence \fref{cnknneooe}. This shows also the close relation between the blow up rate which is far above selfsimilarity\footnote{corresponding to the law $\lambda(t)\sim\sqrt{T-t}$} and the regularity of $u^*$, see \cite{RSz} for related discussions, and explains formally why the problem under consideration should be thought of as "one derivative" above the wave map problem considered in \cite{RR}.\\

{\it 3. Comparison with wave and Schr\"odinger maps}: This result lies in the continuation of the works \cite{RR}, \cite{MRR} on the derivation of stable or codimension one blow up dynamics for the wave map:
\be
\label{wm}
(\mbox{WM})\ \ \left\{\begin{array}{ll} \pa_{tt}u-\Delta u=(|\pa_tu|^2-|\nabla u|^2)u\\ u_{|t=0}=u_0, \ \ \pa_tu_{|t=0}=u_1,\end{array} \right . \ \ (t,x)\in \Bbb R\times \Bbb R^2, \ \ u(t,x)\in \Bbb S^2,
\ee
and the Schr\"odinger map:
\be
\label{nlsmap}
(\mbox{SM})\ \ \left\{\begin{array}{ll} u\wedge\pa_{t}u=\Delta u+|\nabla u|^2u\\ u_{|t=0}=u_0,\end{array} \right . \ \ (t,x)\in \Bbb R\times \Bbb R^2, \ \ u(t,x)\in \Bbb S^2,
\ee
in both cases from $\Bbb R\times \Bbb R^2\to \Bbb S^2$. For the wave map, a stable blow up dynamics within the k-corotational symmetry class \fref{corotshpoere} is exhibited in \cite{RR} for all homotopy number $k\geq 1$ with an almost self similar blow up speed, see also \cite{RS}. In \cite{MRR}, the Schr\"odinger map problem is considered within the k-equivariant symmetry class, ie for solutions of the form 
\be
\label{defucjoei}
v(r,\theta)=e^{k\theta R}w(r),  \ \ w(r)=\left|\begin{array}{lll} w_1(r)\\w_2(r)\\w_3(r)\end{array} \right . 
\ee with 
 \be
 \label{defR}
 R=\left(\begin{array}{lll} 0 & -1 & 0 \\ 1 & 0 &0\\ 0 & 0& 0
\end{array}
\right ),\ \ e^{k\theta R}=\left(\begin{array}{lll} \cos(k\theta) & -\sin(k\theta) & 0 \\ \sin(k\theta)& \cos(k\theta) &0\\ 0 & 0& 1
\end{array}
\right ).
\ee
For $k=1$, a codimension one set of smooth initial data is exhibited for which concentration occurs 
\be
\label{universlbuubling}
\lim_{t\to T}\left\|u(t,r)-e^{\Theta^* R}Q\left(\frac{r}{\lambda(t)}\right)-u^*\right\|_{\mathcal H}=0
\ee 
for some $\Theta^*\in \Bbb R$, $u^*\in \dot{H}^1$ at the speed given by \fref{universallawkgeq}:
$$\lambda(t)=c(u_0)(1+o(1))\frac{T-t}{|\log(T-t)|^2}.$$ Note that the k-equivariant symmetry is also preserved by the wave map and the harmonic heat flow \fref{harmonicheatflow}, and the k-corotational symmetry\footnote{which is not preserved by the Schr\"odinger map} \fref{corotshpoere} corresponds to the k-equivariant symmetry \fref{defucjoei} with $w_2\equiv 0$, and hence such maps are {\it not allowed to rotate around the $e_z$ axis}. Thix extra degree of freedom in k-equivariant symmetry is shown in \cite{MRR} to stabilize the system and leads to a codimension one blow up phenomenon for the Schr\"odinger map. We expect the same phenomenon to occur here, and we conjecture that the blow up solutions constructed in Theorem \ref{thmmain}  correspond to a codimension one phenomenon for the full problem \fref{harmonicheatflow}.\\

{\bf Aknowledgments:} The authors would like to thank Michael Struwe for pointing out to them the relevance of this problem in the continuation of the work \cite{RR}. Part of this work was done while P.R was visiting the ETH, Zurich, which he would like to thank for its kind hospitality. Both authors are supported by the French ERC/ANR project SWAP.\\

{\bf Notations:} We introduce the differential operator $$\Lambda f=y\cdot \nabla f\ \ (\mbox{energy critical scaling}).$$ Given a positive number $b>0$, we let 
\be
\label{defbnot}
B_0=\frac{1}{\sqrt{b}}, \ \ B_1=\frac{|\log b|}{\sqrt{b}}.
\ee
Given a parameter $\lambda>0$, we let $$u_\lambda(r)=u(y)\ \ \mbox{with} \ \ y=\frac{r}{\lambda}.$$  We let $\chi$ be a positive nonincreasing smooth cut off function with $$\chi(y)=\left\{\begin{array}{ll}1\ \ \mbox{for}\ \ y\leq 1,\\ 0\ \ \mbox{for} \ \ y\geq 2.\end{array}\right.$$ Given a parameter $B>0$, we will denote: $$\chi_B(y)=\chi\left(\frac{y}{B}\right).$$
We shall systematically omit the measure in all radial two dimensional integrals and note: $$\int f=\int_0^{+\infty}f(r)rdr.$$


\section{Construction of the approximate profile}


We follow the scheme of proof in \cite{RR}, \cite{MRR} and proceed in this section with the construction of suitable approximate self similar solutions.


\subsection{Asymptotics of the 1-corotational harmonic map}


Let us start with recalling the structure of the harmonic map $Q$ in the context of \fref{mapg} which is the unique -up to scaling- solution to 
\be
\label{harmonicmapequation} 
\Lambda Q=g(Q), \ \ Q(0)=0, \ \ \lim_{r\to +\infty}Q(r)=\pi.
\ee
This equation can be integrated explicitely and leads to the following asymptotics:

\begin{lemma}[Asymptotics of the harmonic map]
\label{harmonicmap}
There holds $Q\in \mathcal C^{\infty}([0,+\infty),[0,\pi))$ with the Taylor expansions\footnote{up to scaling}: 
\be
\label{origin}
Q(y)=\Sigma_{i=0}^pc_{i}y^{2i+1}+O(y^{2p+3})\ \ \mbox{as}\ \ y\to 0,
\ee 
\be
\label{infinity}
Q(y)=\pi-\frac{2}{y}-\Sigma_{i=1}^p\frac{d_i}{y^{2i+1}}+O\left(\frac{1}{y^{2p+3}}\right)\ \ \mbox{as}\ \ y\to +\infty.
\ee
\end{lemma}

\begin{remark} The normalization $$\Lambda Q(y)\sim\frac{2}{y}\ \ \mbox{for}\ \ y\to+\infty$$ from \fref{infinity} has been chosen to match the explicit case of the round sphere: $$Q(y)=2\tan^{-1}(y), \ \ \Lambda Q(y)=\frac{2y}{1+y^2}\sim\frac{2}{y}\ \ \mbox{as}\ \  y\to +\infty.$$
\end{remark}

{\bf Proof of Lemma \ref{harmonicmap}}  From the Taylor expansion of $g$ at $0$ and $\pi$ given by \fref{assumtiong}, we have: $$\int_a^\pi\frac{\pi-\tau-g(\tau)}{(\pi-\tau)g(\tau)}d\tau+\log(\pi-a)\to \left\{\begin{array}{ll} +\infty\ \ \mbox{as}\ \ a\to 0\\ -\infty\ \ \mbox{as}\ \ a\to \pi\end{array}\right.$$ and we may thus find $a\in(0,\pi)$ such that 
\be
\label{normalization}
\int_a^\pi\frac{\pi-\tau-g(\tau)}{(\pi-\tau)g(\tau)}d\tau + \log (\pi - a)=\log 2.
\ee 
We then let
\be
\label{deg} G(u)=\int_b^u\frac{d\tau}{g(\tau)}
\ee
which from \fref{assumtiong} is a diffeomorphism from $(0,\pi)$ onto $(-\infty,+\infty)$. Let now $Q$ be the normalized solution to \fref{harmonicmapequation} given by 
\be
\label{defq}
Q(y)=G^{-1}(\log y), \ \ y\in [0,+\infty).
\ee
We compute near $\pi$ from the normalization \fref{normalization}: 
\bea
\label{developpmntpi}
G(u) & = & \log 2-\log(\pi-u)-\int_u^{\pi}\frac{\pi-\tau-g(\tau)}{(\pi-\tau)g(\tau)}d\tau\\
\nonumber & = & -\log\left(\frac{\pi-u}{2}\right)+\Sigma_{i=0}^p\tilde{d}_i(\pi-u)^{2i+1}+O\left((\pi-u)^{2p+3}\right)\ \ \mbox{as} \ \ y\to \pi
\eea
and near the origin: 
\bee
G(u)& = & \log u+\int_0^u\frac{\tau-g(\tau)}{\tau g(\tau)}d\tau-\int_0^b\frac{\tau-g(\tau)}{\tau g(\tau)}d\tau\\
& = & \log u +\tilde{c}_0+\Sigma_{i=1}^p\tilde{c}_iu^{2i+1}+O\left(u^{2p+3}\right)\ \ \mbox{as}\ \ u\to 0,
\eee
and these developments together with \fref{defq} now yield \fref{origin}, \fref{infinity}. 
This concludes the proof of Lemma \ref{harmonicmap}.


\subsection{The linearized Hamiltonian}


We recall in this section the structure of the linearized operator close to $Q$. Let the potentials 
\be
\label{defpotential}
Z=g'(Q), \ \ V=Z^2+\Lambda Z=f'(Q), \ \ \widetilde{V}=(1+Z)^2-\Lambda Z,
\ee
then the linearized operator close to $Q$ is the Schr\"odinger operator:
\be
\label{defh}
H=-\Delta +\frac{V}{y^2}.
\ee An important consequence of the Bogomolny'i's factorization of the Dirichlet energy \fref{diricheltenergy} is the decomposition  
$$H=A^*A
$$ with $$A= -\partial_y + \frac{Z}{y}, \ \ \ A^*= \partial_y + \frac{1+Z}{y},\ \ Z(y)=g'(Q).$$ The kernels of $A$ and $A^*$ on $\Bbb R^*_+$ are explicit: 
\be
\label{efinitei}
Au=0\ \ \mbox{iff}\ \  u\in \mbox{Span}(\Lambda Q),\ \ A^*u=0\ \ \mbox{iff}\ \ u\in \mbox{Span}\left(\frac{1}{y\Lambda Q}\right),
\ee
and thus the kernel of $H$ on $\Bbb R^*_+$ is: 
\be
\label{kernelh}
Hu=0\ \ \mbox{iff}\ \ u\in \mbox{Span}(\Lambda Q,\Gamma)
\ee 
with 
\be
\label{Gamma} \Gamma(y)=\Lambda Q\int_1^y\frac{dx}{x(\Lambda Q(x))^2}=\left\{\begin{array}{ll} O(\frac1y) \ \ \mbox{as} \ \ y\to 0,\\ \frac{y}{4}+O\left(\frac{\log y}{y}\right)\ \  \mbox{as} \ \ y\to +\infty.\end{array}\right . 
\ee
In particular, $H$ is a positive operator on $\dot{H}^1_{rad}$ with a {\it resonnance} $\Lambda Q$ at the origin induced by the energy critical scaling invariance. We also introduce the conjuguate Hamiltonian
\be
\label{computationhtilde}
\tilde{H}=AA^*=-\Delta +\frac{\widetilde{V}}{y^2}
\ee
which is definite positive by construction and \fref{efinitei}, see Lemma \ref{coerchtilde}. Finally, let us compute using \fref{assumtiong}, \fref{origin}, \fref{infinity} the behavior of $Z, V, \widetilde{V}$ at 0 and $+\infty$ which will be fundamental in our analysis:
\be
\label{comportementz} 
Z(y)=\left\{\begin{array}{ll}1+\Sigma_{i=1}^pc_iy^{2i}+O(y^{2p+2})\ \ \mbox{as}\ \ y\to 0,\\
 				-1+\Sigma_{i=1}^p\frac{c_i}{y^{2i}}+O\left(\frac{1}{y^{2p+2}}\right) \ \ \mbox{as}\ \  y\to+\infty,\end{array}\right.
\ee
\be
\label{comportementv} 
V(y)=\left\{\begin{array}{ll}1+\Sigma_{i=1}^pc_iy^{2i}+O(y^{2p+2})\ \ \mbox{as}\ \ y\to 0,\\
 				1+\Sigma_{i=1}^p\frac{c_i}{y^{2i}}+O\left(\frac{1}{y^{2p+2}}\right) \ \ \mbox{as}\ \  y\to+\infty,\end{array}\right.
\ee
\be
\label{comportementvtilde} 
\widetilde{V}(y)=\left\{\begin{array}{ll}4+\Sigma_{i=1}^pc_iy^{2i}+O(y^{2p+2})\ \ \mbox{as}\ \ y\to 0,\\
 				\Sigma_{i=1}^p\frac{c_i}{y^{2i}}+O\left(\frac{1}{y^{2p+2}}\right) \ \ \mbox{as}\ \  y\to+\infty,\end{array}\right.
\ee
where $(c_i)_{i\geq 1}$ stands for some generic sequence of constants which depend on the Taylor expansion of $g$ at $0$ and $\pi$.

\begin{remark} The exact values for the $\Bbb S^2$ target $g(u)=\sin u$ are given by: $$ Z(y)=\frac{1-y^2}{1+y^2}, \ \ V(y)=\frac{y^4-6y^2+1}{(1+y^2)^2},\ \ \widetilde{V}(y)=2(1+Z)=\frac{4}{1+y^2}.$$
\end{remark}


\subsection{Slowly modulated approximate profiles}


Let $u(t,r)$ be a solution to \fref{mapg}, then the renormalization $$u(t,r)=v(s,y), \ \ y=\frac{r}{\lambda(t)},\ \ \frac{ds}{dt}=\frac{1}{\lambda^2(t)}$$ leads to the self similar equation 
\be
\label{eqselfsimilar}
\pa_sv-\Delta v+b\Lambda v+\frac{f(v)}{y^2}=0, \ \ b=-\lsl.
\ee
We now aim at constructing a suitable approximate solution to \fref{eqselfsimilar} near the harmonic map $Q$ with {\it moderate growth} as $y\to+\infty$ by adapting the slowly modulated ansatz approach developped in \cite{MR2}, \cite{KMR}, \cite{RR}, \cite{MRR}. We will see that this naturally leads to the leading order modulation equation\footnote{see Remark \ref{remarkfon}}:
\be
\label{lawforbzero}
b_s=-b^2(1+o(1)).
\ee

\begin{proposition}[Construction of the approximate profile]
\label{consprofapproch}
Let $M>0$ be a large enough universal constant. Then there exists a small enough universal constant $b^*(M)>0$ such that for all $b \in ]0,b^*(M)[$, there exist profiles $T_1$, $T_2$ et $T_3$, such that
\be
\label{decompqb}
Q_b = Q + bT_1 + b^2T_2 + b^3T_3 = Q + \alpha
\ee
generates an error to the self similar equation \fref{eqselfsimilar} in the regime \fref{lawforbzero}
\be
\label{deferreur}
\Psi_b = - b^2 (T_1 + 2b T_2)  -\Delta Q_b+ b \Lambda Q_b+ \frac{f(Q_b)}{y^2}
\ee
which satisfies:\\
(i) Weighted bounds: 
\be
\label{controleh2erreur}
\int_{y\leq 2B_1}|H\Psi_b|^2 \lesssim b^4|\log b|^2,
\ee
\be
\label{wiehgthtwoewight}
\int_{y\leq 2B_1}\frac{1+|\log y|^2}{1+y^4}|H\Psi_b|^2\lesssim \frac{b^4}{|\log b|^2},
\ee
\be
\label{wiehgthtwoewightthree}
\int_{y\leq 2B_1}\frac{1+|\log y|^2}{1+y^2}|AH\Psi_b|^2\lesssim \frac{b^5}{|\log b|^2},
\ee
\be
\label{controleh4erreur}
\int_{y\leq 2B_1}|H^2\Psi_b|^2 \lesssim \frac{b^6}{|\log b|^2}
\ee
(ii) Flux computation: Let $\Phi_M$ be given by \fref{defdirection}, then:
\be
\label{fluxcomputationone}
\frac{(H\Psi_b,\Phi_M)}{(\Lambda Q, \Phi_M)}=-\frac{2b^2}{|\log b|}+O\left(\frac{b^2}{|\log b|^2}\right),
\ee

\end{proposition}

{\bf Proof of Proposition \ref{consprofapproch}}\\

{\bf step 1} Expansion in powers of $b$.\\

Let us compute the error \fref{deferreur} for a general decomposition \fref{decompqb}. We compute from \fref{decompqb} and a Taylor expansion:
\bee
f(Q_b)& = & f(Q)+b\left[T_1f'(Q)\right]+b^2\left[T_2f'(Q)+\frac12T_1^2f''(Q)\right]\\
& + & b^3\left[T_3f'(Q)+T_1T_2f''(Q)+\frac16f^{(3)}(Q)T_1^3\right] +  R_1+R_2
\eee
with 
\be
\label{defrone}
R_1 =  \frac12f''(Q)\left[\alpha^2-b^2T_1^2-2b^3T_1T_2\right] +  \frac16f^{(3)}(Q)\left[\alpha^3-b^3T_1^3\right],
\ee
\be
\label{defrtwo}
R_2=\frac{\alpha^4}{6}\int_0^1(1-\tau)^3f^{(4)}(Q+\tau\alpha)d\tau.
\ee
Hence from \fref{deferreur}:
\bea
\label{defpsib}
\nonumber \Psi_b 
 &=&b \left( HT_1 + \Lambda Q \right) \\ 
 \nonumber &+& b^2 \left( HT_2 - T_1 + \Lambda T_1 + \frac{f''(Q)}{2y^2} T_1^2  \right) \\ 
 \nonumber &+& b^3\left( HT_3 - 2T_2 + \Lambda T_2 + \frac{f''(Q)}{y^2}T_1T_2 + \frac{1}{6} \frac{f^{(3)}(Q)}{y^2}T_1^3\right) \\ 
&+& b^4\Lambda T_3 + \frac{1}{y^2}\left[R_1+R_2\right].
\eea

{\bf Step 2} Construction of $T_1$.\\

We may invert $H$ explicitely from \fref{kernelh} and a smooth solution at the origin to $Hu=-f$ is given by: 
\be
\label{formulau}
u=\G(y)\int_0^y f\Lambda Qxdx -\Lambda Q(y)\int_{0}^y f\G xdx.
\ee Observe that if $f$ admits the Taylor expansion at the origin $$f(y)=c_1y+c_3y^3+O(y^5),$$ then 
\be
\label{origign}
u(y)=d_3y^3+O(y^5)\ \ \mbox{as} \ \ y\to 0.
\ee Indeed, the Wronskian relation $\Gamma'(\Lambda Q)-(\Lambda Q)'\Gamma=\frac1y$ implies $$A\Gamma=-\Gamma'+\frac{Z}{y}\Gamma=-\Gamma'+\frac{(\Lambda Q)'}{\Lambda Q}\Gamma=-\frac{1}{y\Lambda Q}$$ 
from which using $A\Lambda Q=0$: $$Au=A\Gamma \int_0^y f\Lambda Qxdx=-\frac{1}{y\Lambda Q}\int_0^y f\Lambda Qxdx=c_2y^2+c_4y^4+O(y^6)$$ near the origin for some constants $(c_2,c_4)$. We now integrate using $u=O(y^3)$ at the origin from \fref{formulau} and thus: $$u=-\Lambda Q\int_0^y\frac{Au}{\Lambda Q}dx=d_3y^3+O(y^5).$$
We now let $T_1$ be the solution to $HT_1+\Lambda Q=0$ given by $$ T_{1}(y) = \G(y)\int_0^y (\Lambda Q)^2xdx -\Lambda Q(y)\int_{0}^y \Lambda Q\G xdx.$$ We compute from Lemma \ref{harmonicmap} the behavior:
\be
\label{asymptotetone}
T_1(y)=\left\{\begin{array}{ll} y\log y+e_0y+O\left(\frac{(\log y)^2}{y}\right)\ \ \mbox{as}\ \ y\to+\infty\\ 
						d_3y^3+O(y^5)\ \ \mbox{as}\ \ y\to 0
						\end{array} \right . 
\ee
for some universal constant $e_0$, and similarily:
\be
\label{asympT1}
\Lambda^iT_1(y) = \left\{\begin{array}{ll} y \log y +(e_0+i) y + O\left(\frac{(\log y)^2}{y}\right)\ \ \mbox{as}\ \ y\to+\infty\\ 
						d_{3,i}y^3+O(y^5)\ \ \mbox{as}\ \ y\to 0
						\end{array} \right . \ \ \mbox{for} \ \ 1 \leq i \leq 3 .
\ee

{\bf Step 3} Construction of the radiation $\Sigma_b$.\\ 

Recall the definition \fref{defbnot} and let:
\be
\label{cb}
c_b = \frac{4}{\int \chi_{\frac{B_0}{4}} (\Lambda Q)^2} =\frac{ 2}{ |\log b|}\left( 1 + O( \frac{1}{|\log b|})\right)
\ee
and
\be
\label{db}
d_b = c_b\int_0^{B_0} \chi_{\frac{B_0}{4}} \Lambda Q \G ydy= \frac{C}{b|\log b|}\left( 1 + O\left(\frac{1}{|\log b|}\right)\right)
\ee
Let $\Sigma_b$ be the solution to 
\be
\label{defsigmab}
H\Sigma_b=-c_b\chi_{\frac{B_0}{4}}\Lambda Q +d_b H[(1 - \chi_{3B_0}) \Lambda Q]
\ee
given by 
\be
\label{defsigmabb}
\Sigma_b(y)=\G(y)\int_0^y c_b\chi_{\frac{B_0}{4}}(\Lambda Q)^2xdx -\Lambda Q(y)\int_{0}^y c_{b}\chi_{\frac{B_0}{4}}\G\Lambda Qxdx +d_b (1 - \chi_{3B_0}) \Lambda Q(y)
\ee
Observe that by definition :
\bea
\label{sigmab0}
\Sigma_b = \left\{\begin{array}{ll} c_bT_1\ \ \mbox{for} \ \ y \leq \frac{B_0}{4}\\ \\
4\G \ \ \mbox{for} \ \ y \geq 6B_0.
\end{array} \right. 
\eea
We now estimate for $6B_0 \leq y \leq 2B_1$:
\be
\label{sigmabinfty}
\Sigma_b(y) = y + O(\frac{\log y}{y}) \ \ \ \Lambda \Sigma_b(y) = y + O(\frac{\log y}{y})
\ee
and for $y \leq 6B_0$:
\bea
\label{sigmabmilieu}
\nonumber \Sigma_b(y) & = &c_b\left(\frac{y}{4}+O\left(\frac{\log y}{y}\right)\right)\left[\int_0^y \chi_{\frac{B_0}{4}}(\Lambda Q)^2xdx \right]-c_b\Lambda Q(y)\int_1^{y}O(1)xdx\\
& = & y\frac{\int_0^y\chi_{\frac{B_0}{4}}(\Lambda Q)^2 }{\int \chi_{\frac{B_0}{4}}(\Lambda Q)^2}+O\left(\frac{1+y}{|\log b|}\right)
\eea
and similarily for $y \leq 6B_0$:
\bea
\label{l2sigmabmilieu}
\nonumber \Lambda^i \Sigma_b(y)
& = & y\frac{\int_0^y\chi_{\frac{B_0}{4}}(\Lambda Q)^2 }{\int \chi_{\frac{B_0}{4}}(\Lambda Q)^2}+O\left(\frac{1+y}{|\log b|}\right)\ \ \mbox{for} \ \ 0\leq i\leq 2.
\eea
The equation \fref{defsigmab} and the cancellation $A\Lambda Q=H\Lambda Q=0$ yield the bounds:
\be
\label{boundhigmab}
\int|H\Sigma_b|^2\lesssim \frac{1}{|\log b|},\ \ \int\frac{1+|\log y|^2}{1+y^4}|H\Sigma_b|^2\lesssim \frac{1}{|\log b|^2},
\ee
\be
\label{boundbissigmab}
\int\frac{1+|\log y|^2}{1+y^2}|AH\Sigma_b|^2\lesssim b^2, \ \ \int|H^2\Sigma_b|^2\lesssim \frac{b^2}{|\log b|^2}.
\ee

\begin{bf}
Step 4
\end{bf}
Construction of $T_2$.\\ \\
Let 
\be
\label{defsigma2}
\Sigma_2 = \Sigma_b + T_1 - \Lambda T_1 - \frac{f''(Q)}{2y^2}T_1^2.
\ee
First observe from \fref{assumtiong} that $f=gg'$ is odd and $2\pi$ periodic and thus: $$\forall k \geq 0,  \ \ |f^{(2k)}(u)|+|f^{(2k)}(\pi-u)|\lesssim C_k|u|, \ \ |f^{(2k+1)}(u)|\lesssim 1$$ which implies: 
\be
\label{estimportante}
\forall k\geq 0, \ \ |f^{(2k)}(Q)|\lesssim \frac{y}{1+y^2}.
\ee
We estimate from \fref{Gamma}, \fref{asymptotetone}, \fref{asympT1}, \fref{sigmabinfty}: for $6B_0\leq y\leq 2B_1$, $$\Sigma_2(y)=O\left(\frac{|\log y|^2}{y}\right),$$ and for $y\leq 6B_0$, there holds the behaviour \fref{origign} at the origin and the desired cancellation: 
\bee
\Sigma_2(y)& =&  y\frac{\int_0^y\chi_{\frac{B_0}{4}}(\Lambda Q)^2 }{\int \chi_{\frac{B_0}{4}}(\Lambda Q)^2}-y+O\left(\frac{1+y}{|\log b|}\right)+O\left(\frac{|\log y|^2}{1+y}\right)\\
& = & O\left(\frac{1+y}{|\log b|}(1+|\log (y\sqrt{b})|\right).
\eee
\begin{remark}
\label{remarkfon}
 The above cancellation is due both to the presence of the $T_1$ term in the RHS of \fref{defsigma2} which follows from the choice of modulation equation $b_s=-b^2$ to leading order and cancels the $y\log y$ growth of $T_1$, and the radiation term $\Sigma_b$ which is designed to cancel the remaining $y$ growth in $T_1-\Lambda T_1$.
\end{remark}

We similarily estimate: for $0 \leq i \leq 2$,
\be
\label{sigma2bord}
|\Lambda^i\Sigma_2| \lesssim \frac{y^3}{1+y^2}\left({\bf 1}_{ y \leq 1} + \frac{1 + \log (y\sqrt b)}{|\log b|} {\bf 1}_{1 \leq y \leq 6B_0}\right) + \dfrac{(\log y )^2}{y} {\bf 1}_{y \geq 6B_0}.
\ee
We now let $T_{2}$ be the solution to 
\be
\label{eqtzerotwo}
HT_{2}=\Sigma_{2}
\ee given by 
$$ T_{2}(y)=-\G(y)\int_0^y\Sigma_{2}\Lambda Qxdx +\Lambda Q(y) \int_{0}^y\Sigma_{2}\G xdx
$$
which satisfies \fref{origign} and the estimate from \fref{sigma2bord}:
\be
\label{decayttwo}
\forall y\leq 2B_1, \ \ |\Lambda^iT_{2}(y)|\lesssim  \frac{y^5}{1+y^4}\left({\bf 1}_{ y \leq 1} + \frac1{b|\log b|} {\bf 1}_{y \geq 1}\right),\ \ \mbox{for} \ \ 0 \leq i \leq 3 \\
\ee
We also have the rougher bound: 
\be
\label{roughtwto}
\forall y\leq 2B_1, \ \ |T_2(y)|\lesssim y^3.
\ee

\begin{bf}
Step 5
\end{bf}
Construction of $T_3$.\\
\\
Let 
\be
\label{defsigma3}
\Sigma_3 = 2T_2 - \Lambda T_2 - \frac{f''(Q)}{y^2}T_1T_2 - \frac{1}{6} \frac{f^{(3)}(Q)}{y^2}T_1^3,
\ee
then $\Sigma_3$ satisfies \fref{origign} and we estimate from \fref{asympT1}, \fref{estimportante}, \fref{decayttwo}:
\be
\label{sigma3}
\forall y\leq 2B_1, \ \ |\Lambda^i\Sigma_3(y)|\lesssim  \frac{y^5}{1+y^4}\left({\bf 1}_{ y \leq 1} + \frac1{b|\log b|} {\bf 1}_{y \geq 1}\right),\ \ \mbox{for} \ \ 0 \leq i \leq 2 \\
\ee
We then let $T_3$ be the solution to
 \be
 \label{eqtzerothree}
 HT_{3}=\Sigma_{3}
 \ee given by :
 $$ T_{3}(y)=-\G(y)\int_0^y\Sigma_{3}\Lambda Q+\Lambda Q(y)\int_{0}^y\Sigma_{3}\G
$$
which satisfies \fref{origign} and the estimates from \fref{sigma3}:
\be
\label{T3}
\forall y\leq 2B_1, \ \ | \Lambda^iT_3(y)|\lesssim  \frac{y^7}{1+y^4}\left({\bf 1}_{ y \leq 1} + \frac1{b|\log b|} {\bf 1}_{y \geq 1}\right), \ \ 0\leq i\leq 1,
\ee
\be
\label{roughtwtobis}
|T_3(y)|\lesssim y^3(1+y^2).
\ee
We claim the bounds for $i=0,1$:
\be
\label{estththreone}
\int_{y\leq 2B_1}|H\Lambda^i T_3|^2\lesssim \frac{|\log b|^2}{b^4}, \ \ \int_{y\leq 2B_1}\frac{1+|\log y|^2}{1+y^4}|H\Lambda^i T_3|^2\lesssim \frac{1}{b^4|\log b|^2}
\ee
\be
\label{estththreonebis}
\int\frac{1+|\log y|^2}{1+y^2}|AH\Lambda^i T_3|^2\lesssim \frac{|\log b|^6}{b^2}, \ \ \ \ \int_{y\leq 2B_1}|H^2\Lambda^i T_3|^2\lesssim \frac{1}{b^2|\log b|^2}.
\ee
{\it Proof of \fref{estththreone}, \fref{estththreonebis}}: Observe from  a simple rescaling argument that for any function $f$: 
\be
\label{calculrealcinr}
H\Lambda f=2Hf+\Lambda Hf-\frac{\Lambda V}{y^2}f.
\ee
Hence from \fref{eqtzerothree}, \fref{comportementv}:
$$
 H(\Lambda T_3) = 2 \Sigma_3+ \Lambda \Sigma_3 - \frac{\Lambda V}{y^2} T_3=2 \Sigma_3+ \Lambda \Sigma_3+O(y)
 $$
 \be
 \label{eqhtrop}
 H^2(\Lambda T_3)=H(2\Sigma_3+\Lambda \Sigma_3)+O\left(\frac{1}{1+y}\right).
 \ee
 We thus estimate from \fref{sigma3}, \fref{defbnot} for $i=0,1$:
\be
\label{eniconoefheo} 
\int_{y\leq 2B_1} |H\Lambda^iT_3|^2\lesssim \int_{y\leq 2B_1}\left|\frac{1+y}{b|\log b|} \right|^2\lesssim \frac{B_1^4}{b^2|\log b|^2}\lesssim \frac{|\log b|^2}{b^4},\ee
$$\int_{y\leq 2B_1} \frac{1+|\log y|^2}{1+y^4}|H\Lambda^iT_3|^2\lesssim \int_{y\leq 2B_1}\frac{1}{b^2|\log b|^2}\frac{(1+y^2)(1+|\log y|^2)}{1+y^4}\lesssim \frac{1}{b^4|\log b|^2},
$$
and using the rough bound \fref{roughtwtobis}:
$$\int_{y\leq 2B_1}\frac{1+|\log y|^2}{1+y^2}|AH\Lambda^i T_3|^2\lesssim \int_{y\leq 2B_1} \frac{(1+y^4)(1+|\log y|^2)}{1+y^2}\lesssim B_1^4|\log b|^2\lesssim \frac{|\log b|^6}{b^2}.$$
 The second bound in \fref{estththreonebis} is more subtle and requires further cancellations with respect to \fref{T3}. Indeed, from \fref{calculrealcinr}, \fref{defsigma3}, \fref{eqtzerotwo}:
$$H\Sigma_3=2HT_2-H\Lambda T_2+O\left(\frac{|\log y|^3}{y}\right)=\Lambda \Sigma_2+O\left(\frac{|\log y|^5}{y}\right),$$ $$H\Lambda \Sigma_3=-2\Lambda \Sigma_2-\Lambda^2\Sigma_2+O\left(\frac{|\log y|^5}{y}\right)$$ 
and injecting this into \fref{eqhtrop} with \fref{sigma2bord} yields\footnote{the key here is that $\Sigma_2$ decays at infinity from the cancellation $H\Gamma=0$, and hence the control becomes independent of $B_1$}:
 \bee
 \int_{y\leq 2B_1}|H^2(\Lambda T_3)|^2 & \lesssim &\int_{y\leq 2B_1} \left|\frac{y^3}{1+y^2}\left({\bf 1}_{ y \leq 1} + \frac{1 + \log (y\sqrt b)}{|\log b|} {\bf 1}_{1 \leq y \leq 6B_0}\right) + \dfrac{(\log y )^2}{y} {\bf 1}_{y \geq 6B_0}\right|^2\\
 & \lesssim & \frac{1}{b^2|\log b|^2}
 \eee
and \fref{estththreonebis} is proved.\\
 
 {\bf step 6} Estimate on the error.\\
 
By construction, we have from \fref{defpsib}:
\be
\label{descriptionerreur}
 \Psi_b = b^2 \Sigma_b + b^4\Lambda T_3 + \frac{1}{y^2}(R_1+ R_2).
\ee
We inject into the formulas \fref{defrone}, \fref{defrtwo} the rough bounds \fref{asympT1}, \fref{roughtwto}, \fref{roughtwtobis} and the definition of $B_1$ \fref{defbnot}, and obtain the rough bound: 
$$\forall y\leq 2B_1, \ \ \left|\frac{d^i}{dy^i}\frac{R_1(y)}{y^2}\right|+\left|\frac{d^i}{dy^i}\frac{R_2(y)}{y^2}\right|\lesssim b^4y^{5-i}{\bf 1}_{y\leq 1}+\frac{b^4y^5|\log b|^C}{1+y^{2+i}},  \ \ 0\leq i\leq 4.$$ 
This yields:
$$\int_{y\leq 2B_1}\left|H\left(\frac{R_1}{y^2}\right)\right|^2+\left|H\left(\frac{R_2}{y^2}\right)\right|^2\lesssim b^8|\log b|^C\int_{y\leq 2B_1}\frac{y^{10}}{y^8}\lesssim b^6|\log b|^C,$$
$$\int_{y\leq 2B_1}\left|AH\left(\frac{R_1}{y^2}\right)\right|^2+\left|AH\left(\frac{R_2}{y^2}\right)\right|^2\lesssim b^8|\log b|^C\int_{y\leq 2B_1}\frac{y^{10}}{1+y^{10}}\lesssim b^7|\log b|^C,$$
 $$\int_{y\leq 2B_1} \left|H^2\left(\frac{R_1}{y^2}\right)\right|^2+\left|H^2\left(\frac{R_2}{y^2}\right)\right|^2\lesssim b^8|\log b|^C\int_{y\leq 2B_1}\frac{y^{10}}{1+y^{12}}\lesssim b^8|\log b|^C.$$
Injecting these bounds together with \fref{boundhigmab}, \fref{boundbissigmab}, \fref{estththreone}, \fref{estththreonebis} into \fref{descriptionerreur} yields \fref{controleh2erreur}, \fref{wiehgthtwoewight}, \fref{wiehgthtwoewightthree}, \fref{controleh4erreur}.\\
We now prove the flux computation \fref{fluxcomputationone}. From \fref{descriptionerreur}, \fref{defsigmab}:
 \bee
 \frac{(H\Psi_b,\Phi_M)}{(\Lambda Q, \Phi_M)} &=& \frac{1}{(\Lambda Q, \Phi_M)}\left[\left(-b^2c_b\chi_{\frac{B_0}{4}}\Lambda Q,\Phi_M\right)+O\left(C(M)b^3\right)\right]\\
 & = & -c_bb^2+O\left(C(M)b^3\right)=-\frac{2b^2}{|\log b|}+O\left(\frac{b^2}{|\log b|^2}\right),
\eee
for $|b|<b^*(M)$ small enough. This concludes the proof of Proposition \ref{consprofapproch}


\subsection{Localization of the profile}


We now construct a slowly modulated blow up profile using a localization procedure to avoid artificial unbounded terms for $y\gg B_1$.
\begin{proposition}[Localization of the profile]
\label{localisation}
Let a $\mathcal C^1$ map $ s \mapsto b(s)$ defined on $[0,s_0]$ with a priori bound $\forall s \in [0,s_0]$,
\be
\label{bornes}
0<b(s)<b^*(M),\ \ |b_s| \leq 10 b^2.
\ee
Let the localized profile
$$\tilde{Q}_b(s,y) = Q + b \tilde{T}_1 + b^2\tilde{T}_2 + b^3\tilde{T}_3 = Q + \tilde{\alpha}
$$
where
$$\tilde{T}_i= \chi_{B_1}T_i, \ \ 1\leq i\leq 3.
$$
Then
\be
\label{eqerreurloc}
\partial_s \qbt -\Delta \qbt-\frac {\lambda_s}{\lambda} \Lambda\qbt + \frac{f(\qbt)}{y^2} =\mbox{Mod(t)}+ \tilde{\Psi}_b
\ee
with 
\be
\label{defmode}
Mod(t)=- \left( \frac{\lambda_s}{\lambda} +b\right) \Lambda \qbt + (b_s + b^2)(\tilde{T_1} +2b\tilde{T_2})
\ee
and where $\tilde{\Psi}_b$ satisfies the bounds on $[0,s_0]$ :\\
(i) Weighted bounds: 
\be
\label{controleh2erreurtilde}
\int|H\tilde{\Psi}_b|^2 \lesssim b^4|\log b|^2,
\ee
\be
\label{wiehgthtwoewighttilde}
\int\frac{1+|\log y|^2}{1+y^4}|H\tilde {\Psi}_b|^2\lesssim \frac{b^4}{|\log b|^2},
\ee
\be
\label{wiehgthtwoewightthreetilde}
\int \frac{1+|\log y|^2}{1+y^2}|AH\tilde{\Psi}_b|^2\lesssim \frac{b^5}{|\log b|^2},
\ee
\be
\label{controleh4erreurtilde}
\int|H^2\tilde{\Psi}_b|^2 \lesssim \frac{b^6}{|\log b|^2}
\ee
(ii) Flux computation: Let $\Phi_M$ be given by \fref{defdirection}, then:
\be
\label{fluxcomputationonebis}
\frac{(H\Psit_b,\Phi_M)}{(\Lambda Q, \Phi_M)}=-\frac{2b^2}{|\log b|}+O\left(\frac{b^2}{|\log b|^2}\right),
\ee
\end{proposition}
{\bf Proof of Proposition \ref{localisation}}
\\
\\
{\bf step 1} Localization.
\\
\\
Let $$\Psi^{(1)}_b= - b^2 (\tt_1 + 2b \tt_2)  -\Delta \qbt+ b \Lambda \qbt+ \frac{f(\qbt)}{y^2}.$$ We compute the action of localization which produces an error localized in $[B_1,2B_1]$ up to the term $(1 - \chi_{B_1})\Lambda Q$:
\bea
\label{defpsibloc}
 \nonumber \Psi^{(1)}_b &=& \chi_{B_1} \Psi_b + b(1 - \chi_{B_1})\Lambda Q +  b \Lambda \chi_{B_1}\alpha- 
\alpha \Delta \chi_{B_1} - 2 \partial_y \chi_{B_1} \partial_y\alpha\\
& + & \frac{1}{y^2}\left[f(Q+\chi_{B_1}\alpha)-f(Q)-\chi_{B_1}(f(Q+\alpha)-f(Q))\right].
 \eea
 We estimate in brute force from \fref{asympT1}, \fref{roughtwto}, \fref{roughtwtobis} and the choice of $B_1$: 
 $$\forall y\leq 2B_1, \  \ |\alpha(y)|\lesssim by\left(|\log y|+\frac{by^2}{|\log b|}\right)\lesssim by|\log y|$$ and thus:
 $$\left|b(1 - \chi_{B_1})\Lambda Q +  b \Lambda \chi_{B_1}\alpha- 
\alpha \Delta \chi_{B_1} - 2 \partial_y \chi_{B_1} \partial_y\alpha\right|\lesssim \frac{b}{y}{\bf 1}_{y\geq B_1}+b^2y\log y {\bf 1}_{B_1\leq y\leq 2B_1},
$$
 \bee
 \left| \frac{1}{y^2}\left[f(Q+\chi_{B_1}\alpha)-f(Q)-\chi_{B_1}(f(Q+\alpha)-f(Q))\right]\right| & \lesssim & \frac{|\alpha(y)|}{y^2}|{\bf 1}_{B_1\leq y\leq 2B_1}\\
 & \lesssim & \frac{b\log y}{y}{\bf 1}_{B_1\leq y\leq 2B_1}
 \eee
 from which using \fref{controleh2erreur}, \fref{wiehgthtwoewight}, \fref{wiehgthtwoewightthree}, \fref{controleh4erreur}:
 $$ \int |H\Psi^{(1)}_b|^2\lesssim b^4|\log b|^2+\int_{B_1\leq y\leq 2B_1}\left[\frac{b^2|\log y|^2}{y^{6}}+\frac{b^4|\log y|^2}{y^2}\right]\lesssim b^4|\log b|^2, $$
 $$\int \frac{1+|\log y|^2}{1+y^4}|H\Psi^{(1)}_b|^2\lesssim \frac{b^4}{|\log b|^2}+\int_{B_1\leq y\leq 2B_1}\left[\frac{b^2|\log y|^2}{y^{8}}+\frac{b^4|\log y|^2}{y^4}\right]\lesssim \frac{b^4}{|\log b|^2},$$
 $$\int\frac{1+|\log y|^2}{1+y^2}|AH\Psi^{(1)}_b|^2\lesssim \frac{b^5}{|\log b|^2}+\int_{B_1\leq y\leq 2B_1}\frac{1+|\log y|^2}{1+y^2}\left[\frac{b^2}{y^8}+\frac{b^4|\log y|^2}{y^4}\right]\lesssim \frac{b^5}{|\log b|^2},$$
 \bee
\nonumber   \int|H^2\Psi^{(1)}_b|^2 & \lesssim & \frac{b^6}{|\log b|^2}+\int_{B_1\leq y\leq 2B_1}\left[\frac{b^2|\log y|^2}{y^{10}}+\frac{b^4|\log y|^2}{y^6}\right]\\
 & \lesssim &   \frac{b^6}{|\log b|^2}+\frac{b^2|\log b|^2}{B_1^8}+\frac{b^4|\log b|^2}{B_1^4}\lesssim \frac{b^6}{|\log b|^2}.
 \eee
\begin{remark}
\label{remarkbone}
 This last estimate and \fref{eniconoefheo} govern the choice $B_1=\frac{|\log b|}{\sqrt{b}}$.
\end{remark}
 {\bf step 2} Control of time derivatives.\\
\\
We now compute from \fref{eqerreurloc}:
\be
\label{decomppsib}
\Psit_b = \Psi^{(1)}_b +\tilde{R}, \ \ \tilde{R}= b_s (3b^2 \tilde{T_3} + b \frac{\partial \tilde{T_1}}{\partial b} + b^2 \frac{\partial \tilde{T_2}}{\partial b} + b^3 \frac{\partial \tilde{T_3}}{\partial b})
\ee
and estimate all terms. From \fref{defsigmab}, 
  \be
  \label{cbdbpartialb}
  \frac{\partial c_b}{\partial b}=O\left(\frac{1}{b|\log b|^2}\right), \ \ \frac{\partial d_b}{\partial b}=O\left(\frac{1}{b^2|\log b|}\right)
 \ee
 and thus from \fref{defsigmabb}:
 \bee \frac{\partial \Sigma_{b}}{\partial b}(y) & = & \frac{\partial c_b}{\partial b}T_1{\bf 1}_{y\leq \frac{B_0}{2}}+O\left(\frac{1}{b^2y\log b|}{\bf 1}_{\frac{B_0}{2}\leq y\leq 6B_0}\right)\\
 & = & O\left(\frac{y}{b|\log b|}{\bf 1}_{y\leq B_0}+\frac{1}{b^2y|\log b|}{\bf 1}_{\frac{B_0}{2}\leq y\leq 6B_0}\right).
\eee
We inject this estimate into the explicit formulas for $T_2,T_3$ and conclude:
\be
\label{ttun}
\frac{\partial \tt_1}{\partial b}=  O\left(\frac{y\log y}{b}{\bf 1}_{\frac{B_1}{2}\leq y\leq 2B_1}\right),
\ee
\be
\label{ttdeux}
\frac{\partial \tt_2}{\partial b}=  O\left(\frac{1+y^3}{b|\log b|}{\bf 1}_{y\leq B_0}+\frac{y}{b^2|\log b|}{\bf 1}_{\frac{B_0}{2}\leq y\leq 2B_1}\right),
\ee
\be
\label{tttrois}
\frac{\partial \tt_3}{\partial b}= O\left(\frac{1+y^5}{b|\log b|}{\bf 1}_{y\leq B_0}+\frac{y^3}{b^2|\log b|}{\bf 1}_{\frac{B_0}{2}\leq y\leq 2B_1}\right).
\ee
This yields together with the a priori bound \fref{bornes} the pointwise control:
\bee
|\tilde{R}|\lesssim b^2\left[\frac{b(1+y^3)}{|\log b|}{\bf 1}_{y\leq B_0}+\frac{by^3}{|\log b|}{\bf 1}_{\frac{B_0}{2}\leq y\leq 6B_1}+y\log y{\bf 1}_{\frac{B_1}{2}\leq y\leq 2B_1}\right]
\eee
and hence the bounds:
\bee
\int |H\tilde{R}|^2 & \lesssim &b^4 \int \left|\frac{b(1+y)}{|\log b|}{\bf 1}_{y\leq B_0}+\frac{by}{|\log b|}{\bf 1}_{\frac{B_0}{2}\leq y\leq 6B_1}+\frac{\log y}{y}{\bf 1}_{\frac{B_1}{2}\leq y\leq 2B_1}\right|^2\\
& \lesssim & b^4|\log b|^2,
\eee
\bee
&&\int \frac{1+|\log y|^2}{1+y^4}|H\tilde{R}|^2+ \int\frac{1+|\log y|^2}{1+y^2}|AH\tilde{R}|^2\\
&\lesssim &b^4 \int \frac{1+|\log y|^2}{1+y^4} \left|\frac{b(1+y)}{|\log b|}{\bf 1}_{y\leq B_0}+\frac{by}{|\log b|}{\bf 1}_{\frac{B_0}{2}\leq y\leq 6B_1}+\frac{\log y}{y}{\bf 1}_{\frac{B_1}{2}\leq y\leq 2B_1}\right|^2\\
& \lesssim & \frac{b^5}{|\log b|^2}
\eee
We now track for more cancellation when applying $H^2$. Indeed, from \fref{estththreone}, \fref{T3}:
$$ \int|H^2\tt_3|^2\lesssim \frac{1}{b^2|\log b|^2}+\int\left|\frac{y^3}{by^4|\log b|}{\bf 1}_{B_1\leq y\leq 2B_1}\right|^2\lesssim \frac{1}{b^2|\log b|^2}.$$ From direct inspection:
$$\int\left| H^2\left(\frac{\partial T_1}{\partial b}\right)\right|^2\lesssim \int_{B_1\leq y\leq 2B_1}\left|\frac{\log y}{by^3}\right|^2\lesssim \frac{|\log b|^2}{b^2B_1^4}\lesssim \frac{1}{|\log b|^2}.$$ Next:
$$H^2\left(\frac{\partial \tt_2}{\partial b}\right)=\chi_{B_1}H^2\left(\frac{\partial T_2}{\partial b}\right)+O\left(\frac{1}{y^3b^2|\log b|}{\bf 1}_{\frac{B_1}{2}\leq y\leq 2B_1}\right)$$
and from \fref{eqtzerotwo}, \fref{defsigma2}:
\bee
H^2\left(\frac{\partial T_2}{\partial b}\right) & = & H\left(\frac{\partial \Sigma_2}{\partial b}\right)=H\left(\frac{\partial \Sigma}{\partial b}\right)\\
& = & O\left(\frac{1}{b(1+y)|\log b|^2}{\bf 1}_{y\leq 2B_0}+\frac{1}{b(1+y)|\log b|}{\bf 1}_{\frac{B_0}{2}\leq y\leq 2B_0}\right)
\eee
which yields the bound:
\bee
\int\left|H^2\left(\frac{\partial \tt_2}{\partial b}\right) \right|^2& \lesssim & \int \left|\frac{1}{b(1+y)|\log b|^2}{\bf 1}_{y\leq 2B_0}+\frac{1}{b(1+y)|\log b|}{\bf 1}_{\frac{B_0}{2}\leq y\leq 2B_0}\right.\\
& + & \left .\frac{1}{y^3b^2|\log b|}{\bf 1}_{\frac{B_1}{2}\leq y\leq 2B_1}\right|^2\lesssim \frac{1}{b^2|\log b|^2}.
\eee
Similarily, 
\bee
H^2\left(\frac{\partial \tt_3}{\partial b}\right)& = & \chi_{B_1}H^2\left(\frac{\partial T_3}{\partial b}\right)+O\left(\frac{1}{yb^2|\log b|}{\bf 1}_{\frac{B_1}{2}\leq y\leq 2B_1}\right)\\
& = & O\left(\frac{1+y}{b|\log b|}{\bf 1}_{y\leq B_0}+\frac{1}{yb^2|\log b|}{\bf 1}_{\frac{B_0}{2}\leq y\leq 2B_0}+\frac{1}{yb^2|\log b|}{\bf 1}_{\frac{B_1}{2}\leq y\leq 2B_1}\right)
\eee
from which:
\bee
\int\left|H^2\left(\frac{\partial \tt_3}{\partial b}\right) \right|^2 & \lesssim & \int\left|\frac{1+y}{b|\log b|}{\bf 1}_{y\leq B_0}+\frac{1}{yb^2|\log b|}{\bf 1}_{\frac{B_0}{2}\leq y\leq 2B_0}+\frac{1}{yb^2|\log b|}{\bf 1}_{\frac{B_1}{2}\leq y\leq 2B_1}\right|^2\\
& \lesssim & \frac{1}{b^4|\log b|^2}.
\eee
We inject these estimates into \fref{decomppsib} and obtain: $$\int |H^2\tilde{R}|^2\lesssim \frac{b^6}{|\log b|^2}.$$Injecting the collection of estimates of step 1 and step 2 into \fref{decomppsib} now yields the expected bounds \fref{controleh2erreurtilde}, \fref{wiehgthtwoewighttilde}, \fref{wiehgthtwoewightthreetilde}, \fref{controleh4erreurtilde}.\\

 {\bf step 3} Flux computation.\\
 
 By definition, $\Psi^{(1)}_b=\Psi_b$ on $\mbox{Supp}(\Phi_M)\subset[0,2M],$ and from \fref{ttun}, \fref{ttdeux},  \fref{tttrois} and \fref{T3}: $$|\tilde R|\lesssim C(M)|b_s|b\lesssim C(M)b^3\ \ \mbox{on Supp}(\Phi_M).$$ This estimate together with \fref{fluxcomputationone} now yields \fref{fluxcomputationonebis}.\\
This concludes the proof of Proposition \ref{localisation}.\\

We introduce a second localization of the profile near $B_0$ which will be used only to capture some further cancellation in the proof of the bound \fref{regularityustar}.

\begin{lemma}[Second localization]
\label{localisationbis}
Let a $\mathcal C^1$ map $ s \mapsto b(s)$ defined on $[0,s_0]$ with a priori bound \fref{regularityustar}. Let the localized profile
\be
\label{defhatqb}
\hat{Q}_b(s,y) = Q + b \hat{T}_1 + b^2\hat{T}_2 + b^3\hat{T}_3 = Q + \hat{\alpha}
\ee
where
$$\hat{T}_i= \chi_{B_0}T_i, \ \ 1\leq i\leq 3.
$$
Let the radiation:
\be
\label{defraidiation}
\zeta_b=\tilde{\alpha}-\hat{\alpha}
\ee
and the error
$$
\partial_s \qbh -\Delta \qbh-\frac {\lambda_s}{\lambda} \Lambda\qbh + \frac{f(\qbh)}{y^2} =\widehat{Mod}(t)+ \hat{\Psi}_b
$$
with 
\be
\label{defmodebis}
\widehat{Mod}(t)=- \left( \frac{\lambda_s}{\lambda} +b\right) \Lambda \qbh + (b_s + b^2)(\hat{T_1} +2b\hat{T_2}).
\ee
Then there holds the bounds:
\be
\label{radiation}
\int|H\zeta_b|^2\lesssim b^2|\log b|^2,\ \ \Sigma_{i=0}^2\int\frac{|\pa_y^i\zeta_b|^2}{1+y^{8-2i}}\lesssim b^4|\log b|^C,
\ee
\be
\label{radiationbis}
\Sigma_{i=0}^2\int\frac{|\pa_y^i\zeta_b|^2}{1+y^{2(3-i)}}\lesssim b^3|\log b|^C,
\ee
\be
\label{controleh2erreurtildehat}
Supp(H\hat{\Psi}_b)\subset [0,2B_0] \ \ \mbox{and} \ \ \int|H\hat{\Psi}_b|^2 \lesssim b^4|\log b|^2.
\ee
\end{lemma}

\begin{remark} Note that this localization near $B_0$ displays the same properties like the one near $B_0$ at the $\mathcal H^2$ level and \fref{controleh2erreurtilde}, \fref{controleh2erreurtildehat} are comparable. The estimate corresponding to \fref{controleh4erreurtilde} would however be worse for $\hat{\Psi}_b$ due to the terms induced by localization, see Remark \ref{remarkbone}. Hence we will use the $B_1$ localization to control high derivatives norms, see Proposition \ref{AEI2} and the control of $\mathcal E_4$ below, and $B_0$ localization for lower order control, see section \ref{sectionbootstrap} and the control of $\mathcal E_2$.
\end{remark}

{\bf Proof of Lemma \ref{localisationbis}} By construction, $$\zeta_b=(\chi_{B_1}-\chi_{B_0})(bT_1+b^2T_2+b^3T_3)$$ and thus from \fref{decayttwo}, \fref{T3}: $$\int|H\zeta_b|^2\lesssim \int_{y\leq 2B_1}\left|\frac{by\log y}{y^2}+\frac{b^2y}{by^2|\log b|}+\frac{b^3y^3}{by^2|\log b|}\right|^2\lesssim b^2|\log b|^2,$$
$$\Sigma_{i=0}^2\int\frac{|\pa^i_y\zeta_b|^2}{1+y^{8-2i}}\lesssim \int_{B_0\leq y\leq 2B_1}\frac{1}{1+y^8}\left[b^2y^2|\log y|^2+\frac{b^6y^{6}}{b^2|\log b|^2}\right]\lesssim b^4|\log b|^C,$$
$$\Sigma_{i=0}^2\int\frac{|\pa_y^i\zeta_b|^2}{1+y^{2(3-i)}}\lesssim \int_{B_0\leq y\leq 2B_1}\frac{1}{1+y^6}\left[b^2y^2|\log y|^2+\frac{b^6y^{6}}{b^2|\log b|^2}\right]\lesssim b^3|\log b|^C.$$
 The localization property \fref{controleh2erreurtildehat} directly follows from the analogue of the formula \fref{defpsibloc}, \fref{decomppsib} for $\hat{\Psi}_b$ and the cancellation $H\Lambda Q=0$, while the proof of the estimate \fref{controleh2erreurtildehat} is very similar to the one of \fref{controleh2erreurtilde} and left to the reader. 

This concludes the proof of Lemma \ref{localisationbis}.


  \section{The trapped regime}


This section is devoted to the description of the set of initial data and the corresponding trapped regime in which the singularity formation described by Theorem \ref{thmmain} will occur.


\subsection{Setting the bootstrap}


We describe in this section the set of initial data leading to the blow up scenario of Theorem  \ref{thmmain}. Let a 1-corotational map $$v_0\in \dot{H}^1\cap \dot{H}^4\ \ \mbox{with} \ \ \|\nabla (v_0-\mathcal Q)\|_{L^2}\ll 1$$ and $v(t)\in \mathcal C([0,T), \dot{H}^1\cap \dot{H}^4)$, $0<T\leq +\infty$ be the corresponding solution to \fref{hfgeneral}. First recall from standard argument the blow up criterion:
\be
\label{regularitiyt}
T<+\infty\ \ \mbox{implies}\ \ \|\Delta v(t)\|_{L^2}\to +\infty\ \ \mbox{as} \ \ t\to T.
\ee
From Lemma \ref{regularitycorot}, $v$ admits on a small time interval $[0,t_1]$ a decomposition 
\be
\label{defvvu}
v(t,x)=\left|\begin{array}{lll} g(u(t,r))\cos\theta\\g(u(t,r))\sin\theta\\z(u(t,r))\end{array}\right.
\ee where $\tilde{\e}(t,r)=u(t,r)-Q(r)$ satisfies the boundary condition \fref{boudnaryconditions} and the regularity \fref{sobolefofiue}, \fref{sobolefofiue} displayed in Lemma \ref{regularitycorot}. Moreover, from the initial smallness $|\nabla \et(0)\|_{L^2}\ll 1$, we may from standard modulation argument introduce the unique decomposition 
\be
\label{decompu}
u(t,r)=(\tilde{Q}_{b(t)}+\e(t,r))_{\lambda(t)}, \ \ \|\nabla \e(t)\|_{L^2}+|b_(t)|\ll 1, \ \ \lambda(t)>0
\ee where $\e(t)$ satisfies the orthogonality condition: 
\be
\label{ortho}
\forall t\in [0,t_1], \ \ (\e(t),\Phi_M)=(\e(t),H\Phi_M)=0.
\ee
Here given $M>0$ large enough, $\Phi_M$ corresponds to the fixed direction
\be
  \label{defdirection}
  \Phi_M = \chi_M \Lambda Q - c_M H(\chi_M\Lambda Q)
    \ee
  with
 $$c_M = \dfrac{\left(\chi_M \Lambda Q,T_1\right)}{\left(H\left( \chi_M \Lambda Q \right), T_1\right)} = c_{\chi} \frac{M^2}{4}(1+o_{M \to + \infty}(1)).$$
Observe by construction that 
 \be
 \label{intM}
 \int |\Phi_M|^2 \lesssim |\log M|,\ \ (\Phi_M,T_1)=0
 \ee
and the scalar products
\be
\label{prodscalar}
(\Lambda Q, \Phi_M) = (-HT_1,\Phi_M) = (\chi_M\Lambda Q,\Lambda Q) = 4 \log M (1+o_{M \to + \infty}(1))
\ee
are non degenerate. The existence of the decomposition \fref{decompu} is then a standard consequence of the implicit function theorem and the explicit relations $$\left(\frac{\pa}{\pa\l}(\qbt)_\lambda,\ \frac{\pa}{\pa b}(\qbt)_\lambda\right)|_{\l=1,b=0}=(\Lambda Q,T_1)$$ which ensure the nondegeneracy of the Jacobian:
\bee
\left|\begin{array}{ll} (\frac{\pa}{\pa\l}(\qbt)_\lambda,\Phi_M)& (\frac{\pa}{\pa b}(\qbt)_\lambda,\Phi_M)\\  (\frac{\pa}{\pa\l}(\qbt)_\lambda,H\Phi_M)& (\frac{\pa}{\pa b}(\qbt)_\lambda,H\Phi_M)\end{array}\right |_{\l=1,b=0}& = & \left|\begin{array}{ll} (\Lambda Q,\Phi_M)& 0 \\ 0 & (T_1,H\Phi_M)\end{array}\right |\\
& = & -(\Lambda Q,\Phi_M)^2\neq 0.
\eee 
From Lemma \ref{regularitycorot}, we may measure the regularity of the map through the following norms of $\e$: the energy norm
\be
\label{defenergyspcar}
\|\e\|_{\mathcal H}^2=\int|\pa_y\e|^2+\int\frac{|\e|^2}{y^2},
\ee
and higher order Sobolev norms adapated to the linearized operator 
\be
\label{defnorme}
\mathcal E_{2k}=\int|H^k\e|^2, \ \ 1\leq k\leq 2.
\ee
We now assume the following bounds on initial data which describe an open $\dot{H}^1\cap H^4$ affine space of 1-corotational initial data around $\mathcal Q$:
\begin{itemize}
\item Smallness and positivity of b(0) :
\be
\label{init1}
0 <b(0)<b^*(M)\ll1.
\ee
\item Smallness of the excess of energy:
\be
\label{init2energy}
\int | \nabla \varepsilon (0)|^2 + \int \left|\dfrac{ \varepsilon (0)}{y}\right|^2 \leq b^2(0),
\ee
\be
\label{init2}
|\mathcal E_2(0)| + |\mathcal E_4(0)| \leq [b(0)]^{10}.
\ee
\end{itemize}
The propagation of the $\dot{H}^1\cap \dot{H}^4$ regularity by the parabolic heat flow and Lemma  \ref{regularitycorot} ensure that these estimates hold on some small enough time interval $[0,t_1]$. From standard argument, there also holds the regularity $(\lambda, b)\in \mathcal C^1([0,t_1],\Bbb R^*_+\times \Bbb R)$. Given a large enough universal constant $K>0$ -independent of $M$-, we assume on $[0,t_1]$ the following bootstrap bounds:
\begin{itemize}
\item Control of b(t) :
\be
\label{init1h}
0 <b(t)<10b(0).
\ee
\item Control of the radiation:
\be
\label{init2h}
\int | \nabla \varepsilon (t)|^2 + \int \left|\dfrac{ \varepsilon (t)}{y}\right|^2 \leq 10\sqrt{b(0)},
\ee
\be
\label{init3h}
|\mathcal E_2(t)| \leq K b^2(t)|\log b(t)|^5,
\ee
\be
\label{init3hbis}
|\mathcal E_4(t)| \leq K \frac{b^4(t)}{|\log b(t)|^2}.
\ee
\end{itemize} 
The following proposition describes the contraction of the bootstrap regime and is the core of the proof.

\begin{proposition}[Bootstrap control of $b$ and $\varepsilon$]
\label{bootstrap}
Assume that K in \fref{init1h}, \fref{init2h}, \fref{init3h} and \fref{init3hbis} has been chosen large enough. Then, $\forall t \in [0,t_1]$ :
\be
\label{init1hb}
0 <b(t)<2b(0),
\ee
\be
\label{init2hb}
\int | \nabla \varepsilon (t)|^2 + \int \left|\dfrac{ \varepsilon (t)}{y}\right|^2 \leq \sqrt{b(0)},
\ee
\be
\label{init3hb}
|\mathcal E_2(t)| \leq \frac K2 b^2(t)|\log b(t)|^5,
\ee
\be
\label{init3hbisb}
|\mathcal E_4(t)| \leq \frac K2 \frac{b^4(t)}{|\log b(t)|^2}.
\ee
\end{proposition}

The rest of this section is devoted to developping the tools needed for the proof of Proposition \ref{bootstrap} which will be completed in section \ref{sectionbootstrap}.


\subsection{Equation for the radiation}


Recall the decomposition of the flow:
  \bee
  u(t,r) = (\tilde{Q}_{b(t)} + \varepsilon) (t,\frac{r}{\lambda(t)}) = (Q + \tilde{\alpha}_{b(t)})_{\lambda(t)} + w(t,r).
  \eee
 We introduce the rescaled time
  \bea
  \nonumber s(t) = \int_0^t \frac{d \tau}{\lambda^2(\tau)}
  \eea
  and use the rescaling formulas
  \bea
  \nonumber u(t,r) = v(s,y)\mbox{,} \ \ \ y = \frac{r}{\lambda(t)} \mbox{,}  \ \ \ \ \ \partial_t u = \frac{1}{\lambda^2(t)}(\partial_s v - \frac{\lambda_s}{\lambda}\Lambda v)_{\lambda}
  \eea  
  to derive the equation for $\e$ in renormalized variables:
  \be
  \label{eqepsilon}
  \partial_s \varepsilon  - \frac{\lambda_s}{\lambda} \Lambda \varepsilon + H \varepsilon = F - Mod = \mathcal F.
  \ee
  Here  $H$ is the linearized operator given by \fref{defh}, $\mbox{Mod}(t)$ is given by \fref{defmode},
  \be
  \label{defF}
  F = -\Psit_b +L(\varepsilon) - N(\varepsilon)
  \ee
where $L$ is the linear operator corresponding to the error in the linearized operator from $Q$ to $\tilde{Q}_b$:
  \be
  \label{defLe}
 L(\varepsilon) =  \frac{f'(Q)-f'(\qbt)}{y^2}\e
  \ee
  and the remainder term is the purely nonlinear term:
  \be
  \label{defNe}
  N(\varepsilon) = \frac{f(\tilde{Q}_b+ \varepsilon) - f(\tilde{Q}_b)-\e f'(\qbt)}{y^2}. 
  \ee
We also need to write the flow \fref{eqepsilon} in original variables. For this, let the rescaled operators
$$A_{\lambda} = -\partial_r + \frac{Z_{\lambda}}{r}, \ \ \ A^*_{\lambda} = \partial_r + \frac{1+Z_{\lambda}}{r}\ \ $$
\be
\label{defhtilde}
H_{\lambda} = A^*_{\lambda}A_{\lambda} = -\Delta + \frac{V_{\lambda}}{r^2}, \ \ \tilde{H}_{\lambda} = A_{\lambda}A^*_{\lambda} = -\Delta + \frac{\tilde{V}_{\lambda}}{r^2},
\ee
and the renormalized function $$w(t,r)=\e(s,y),$$
then \fref{eqepsilon} becomes:
  \be
  \label{eqenwini}
  \partial_t w + H_{\lambda}w = \frac 1{\lambda^2} \mathcal F_{\lambda}
  \ee
  

    \subsection{Modulation equations}
  

Let us now compute the modulation equations for $(b,\l)$ as a consequence of the choice of orthogonality conditions \fref{ortho}.

  \begin{lemma}[Modulation equations]
\label{modulationequations}
There holds the bound on the modulation parameters :
\be
\label{parameters}
\left|\frac{\lambda_s}{\lambda} + b\right| \lesssim \frac{b^2}{|\log b|} + \frac{1}{\sqrt{\log M}}\sqrt{\mathcal E_4},
\ee
\be
\label{parameterspresicely}
\left| b_s + b^2\left(1 + \frac{2}{|\log b|}\right) \right| \lesssim \frac{1}{\sqrt{\log M}} \left( \sqrt{ \mathcal E_4} +  \frac{b^2}{|\log b|}  \right).
\ee
\end{lemma}

\begin{remark} Note that this implies in the bootstrap the rough bounds:
\be
\label{rougboundpope}
|b_s|+\left|\lsl+b\right|\leq 2b^2.
\ee
and in particular \fref{bornes} holds. 
\end{remark}

{\bf Proof of Lemma \ref{modulationequations}}\\

{\bf step 1} Law for b.\\

Let
$$V(t) = |b_s + b^2| + \left|\frac{\lambda_s}{\lambda} + b\right|$$
We take the inner product of \fref{eqepsilon} with $H\Phi_M$ and compute:
\bea
\label{modequone}
(\mbox{Mod}(t),H\Phi_M)& = & -(\Psit_b,H\Phi_M)-(H\e,H\Phi_M)\\
\nonumber & - & \left(-\lsl\Lambda \e-L(\e)+N(\e),H\Phi_M\right).
\eea
 We first compute from the construction of the profile, \fref{defmode} and the localization $\mbox{Supp}(\Phi_M)\subset[0,2M]$ from \fref{defdirection}:
 \bee
\nonumber \left(H\left(Mod(t)\right),\Phi_M\right) &=& -\left(b+\frac{\lambda_s}{\lambda}\right) \left(H\Lambda \qbt , \Phi_M\right)+ \left(b_s + b^2\right)\left(H\left(\tilde{T}_1 + 2b \tilde{T}_2\right),\Phi_M\right)\\
& = &  -(\Lambda Q, \Phi_M)(b_s + b^2)+O\left(c(M)b| V(t)|\right).
\eee
The linear term in \fref{modequone} is estimated\footnote{Observe that we do not use the interpolated bounds of Lemma \ref{lemmainterpolation} but directly the definition \fref{defnorme} of $\mathcal E_4$, and hence the dependence of the constant in $M$ is explicit what is crucial for the analysis.} from \fref{intM}:
$$
\nonumber \left| (H \varepsilon,H\Phi_M) \right| \lesssim \| H^2\varepsilon \|_{L^2}\sqrt{\log M}=\sqrt{\log M\mathcal E_4}
$$
and the remaining nonlinear term is estimated using the Hardy bounds of Appendix A:
$$\left|\left(-\lsl\Lambda \e+L(\e)+N(\e),H\Phi_M\right)\right|\lesssim C(M)b(\sqrt{\mathcal E_4}+|V(t)|).$$
We inject these estimates into \fref{modequone} and conclude from \fref{prodscalar} and the fundamental flux computation \fref{fluxcomputationonebis}:
\bee
b_s+b^2& = & \frac{(\Psit_b,H\Phi_M)}{(\Lambda Q,\Phi_M)}+O\left(\frac{\sqrt{\log M\mathcal E_4}}{\log M}\right)+O\left(C(M)b|V(t)\right)
\eee
and hence the first modulation equation:
\bea
\label{firstmodulationequation} b_s + b^2 = -\frac{2b^2}{|\log b|}\left( 1 + O\left(\dfrac{1}{|\log b|}\right) \right) + O\left( \sqrt{\dfrac{\mathcal E_4}{\log M}}+C(M)b|V(t)|\right).
\eea

{\bf step 2} Degeneracy of the law for $\lambda$.\\

We now take the inner product of \fref{eqepsilon} with $\Phi_M$ and obtain:
$$(\mbox{Mod}(t),\Phi_M)=-(\Psit_b,\Phi_M)-(H\e,\Phi_M)-\left(-\lsl\Lambda \e+L(\e)+N(\e),\Phi_M\right).$$
Note first that the choice of orthogonality conditions \fref{ortho} gets rid of the linear term in $\varepsilon$: $$(H\e,\Phi_M)=0.$$ Next, we compute from  \fref{prodscalar} and the orthogonality \fref{intM}:
\bee
\nonumber \left(\mbox{Mod}(t),\Phi_M\right) &=& -\left(b+\frac{\lambda_s}{\lambda}\right) \left(\Lambda \qbt , \Phi_M\right)+ \left(b_s + b^2\right)\left(\tilde{T}_1 + 2b \tilde{T}_2,\Phi_M\right)\\
& = &  - 4 \log M (1+o_{M \to + \infty}(1))\left(\lsl+b\right)+O\left(C(M)b| V(t)|\right).
\eee
and observe the cancellation from \fref{sigmab0}, \fref{intM}:
$$
\nonumber \left| \left(\Psit_b,\Phi_M \right) \right| \lesssim b^2|(\Sigma_b, \Phi_M)| + O(C(M)b^3)= c_b b^2 |(T_1,\Phi_M)| +O(C(M)b^3)=O(C(M)b^3).$$
Nonlinear terms are easily estimated using the Hardy bounds of Appendix A:
$$\left|\left(-\lsl\Lambda \e+L(\e)+N(\e),\Phi_M\right)\right|\lesssim C(M)b\left(\sqrt{\mathcal E_4}+|V(t)|\right).$$

We thus obtain the modulation equation for scaling:
$$\left|\lsl+b\right|\lesssim b^3C(M)+bC(M)O\left(\sqrt{\mathcal E_4}+ |V(t)|\right).
$$
Combining this with \fref{firstmodulationequation}  yields the bound $$|V(t)|\lesssim \frac{b^2}{|\log b|} + \frac{1}{\sqrt{log M}}\sqrt{\mathcal E_4}$$ which together with \fref{firstmodulationequation} again now implies the refined bound \fref{parameterspresicely}. This concludes the proof of Lemma \ref{modulationequations}.


\subsection{The Lyapounov monotonicity}


We now turn to the core of the argument which is the derivation of a suitable Lyapounov functional at the $H^4$ level. The parabolic structure will yield further dissipation with respect to the analysis of dispersive problems in \cite{RR}, \cite{MRR}, what will allow us to treat a general metric $g$.

\begin{proposition}[Lyapounov monotonicity]
\label{AEI2}
There holds:
\bea
\label{monoenoiencle}
\frac{d}{dt} \left\{\frac{1}{\lambda^6}\left[\mathcal E_4+O\left(\sqrt b\frac{b^4}{|\log b|^2}\right)\right]\right\} \leq C\frac b {\lambda^8}\left[ \frac{\mathcal E_4}{\sqrt{\log M}}+\frac{b^4}{|\log b|^2}+\frac{b^2}{|\log b|}\sqrt{\mathcal E_4} \right]
\eea
for some universal constant $C>0$ independent of $M$ and of the bootstrap constant $K$ in \fref{init1h}, \fref{init2h}, \fref{init3h}, \fref{init3hbis}, provided $b^*(M)$ in \fref{init1} has been chosen small enough.
\end{proposition}

{\bf Proof of Proposition \ref{AEI2}}\\

The proof relies on the derivation of the energy identity for suitable derivatives of $\e$ seen in original variables ie $w$, and repulsivity properties of the corresponding time dependent Hamiltonian $\tilde{H}_\l$. The control of the solution is then ensured thanks to coercivity properties of the iterated Hamiltonian $H,H^2$ under the orthogonality conditions \fref{ortho}, see Lemma \ref{lemmacoer}. Nonlinear terms will be estimated using the interpolated bounds of Lemma \ref{lemmainterpolation} which will be implicitely used all along the proof.\\

{\bf step 1} Suitable derivatives.\\

We define the derivatives of $w$ associated with the linearized Hamiltonian $H$:
$$w_1 = A_{\lambda}w, \ \ \ w_2 = A^*_{\lambda}w_1, \ \ \ w_3 = A_{\lambda}w_2 
$$
which satisfy from \fref{eqenwini}:
$$\partial_{t} w_1 + \tilde{H}_{\lambda} w_1 = \frac{\partial_t Z_\l}{r}w + A_{\lambda} \left(\frac 1{\lambda^2} \mathcal F_{\lambda}\right)
$$
\be
\label{eqw2}
\partial_{t} w_2 + H_{\lambda} w_2 = \frac{\partial_{t} V_{\lambda}}{r^2} w + H_{\lambda} \left(\frac 1{\lambda^2} \mathcal F_{\lambda}\right)
\ee
\be
\label{eqw3}
\partial_{t} w_3 + \tilde{H}_{\lambda} w_3 = \frac{\partial_{t} Z_{\lambda}}{r} w_2 + A_{\lambda} \left( \frac{\partial_{t} V_{\lambda}}{r^2} w \right) + \left(AH \right)_{\lambda} \left(\frac 1{\lambda^2} \mathcal F_{\lambda}\right)
\ee
We similarily use in the following steps the notation: $$\e_1=A\e, \ \ \e_2=A^*\e_1, \ \ \e_3=A\e_2.$$ Observe from \fref{defnorme} that
\be
\label{icneoneohie}
\mathcal E_4=\int|A^*\e_3|^2.
\ee
We recall the action of time derivatives on rescaling: $$\partial_tv_{\lambda}=\frac{1}{\l^2}\left(\pa_sv-\lsl\Lambda v\right)_{\lambda}.$$

{\bf step 2} Modified energy identity.\\

We compute the energy identity on \fref{eqw3} using \fref{defhtilde}:
\par
\bea
\label{firstestimate}
\nonumber  \frac{1}{2}\frac{d \mathcal E_4}{dt} & = & \frac{1}{2} \frac{d}{d t}\left\{ \int \tilde{H}_{\lambda}w_3 w_3 \right\}= \int \tilde{H}_{\lambda} w_3 \partial_{t} w_3 + \int \frac{\partial_{t} \tilde{V}_\l}{2r^2} w_3^2\\
\nonumber & = & - \int (\tilde{H}_\l w_3)^2+b\int \frac{(\Lambda \tilde{V})_\l}{2\l^2r^2} w_3^2-\left(\lsl+b\right)\int \frac{(\Lambda \tilde{V})_\l}{2\l^2r^2} w_3^2 \\
& + &\int \tilde{H}_\l w_3\left[\frac{\partial_{t} Z_{\lambda}}{r} w_2 +  A_{\lambda}\left(\frac{\partial_{t} V_{\lambda}}{r^2}w \right)+\left(AH \right)_{\lambda}\left( \frac{1}{\lambda^2}\mathcal F_{\lambda}\right)\right].
\eea
We now aim at using the dissipative term $\int(\tilde{H}w_3)^2$ to treat the quadratic terms in the RHS of \fref{firstestimate}. Observe however that this quantity is delicate to use because it is positive but not coercive a priori\footnote{as can be seen by considering the zero of $\tilde{H}$ given $w_3=\frac{1}{y\Lambda Q}\int_0^y\tau(\Lambda Q)^2d\tau$.}. Also one can explicitely compute for the sphere target
$$b\frac{(\Lambda V)_\l}{\l^2r^2}=-\frac{b}{\l^8}\frac{4}{(1+y^2)^2}<0,$$  
and thus the critical in size quadratic term has the right sign in this case: $$b\int \frac{(\Lambda \tilde{V})_\l}{\l^2r^2} w_3^2<0$$ which would allow some simplification of our analysis. However, this sign property does not seem to hold a priori for the general metric $g$ we consider. We nevertheless claim using a similar algebra as in \cite{RR}, \cite{MRR} that this term can be treated thanks to a further integration by parts in time which in the dispersive cases would correspond to a Morawetz type computation. Indeed, we compute from \fref{eqw2}, \fref{eqw3}:
\bee
\frac{d}{dt}\left\{\int \frac{b(\Lambda Z)_\l}{\l^2r}w_3w_2\right\}& = & \int \frac{d}{dt}\left(\frac{b(\Lambda Z)_\l}{\l^2r}\right)w_3w_2\\
& + & \int \frac{b(\Lambda Z)_\l}{\l^2r}w_2\left[-\tilde{H}_{\l}w_3+\frac{\partial_{t} Z_{\lambda}}{r} w_2 + A_{\lambda} \left( \frac{\partial_{t} V_{\lambda}}{r^2} w \right) + \left(AH \right)_{\lambda} \left(\frac 1{\lambda^2} \mathcal F_{\lambda}\right)\right]\\
& + & \int\frac{b(\Lambda Z)_\l}{\l^2r}w_3\left[-A_\l^*w_3+\frac{\partial_{t} V_{\lambda}}{r^2} w + H_{\lambda} \left(\frac 1{\lambda^2} \mathcal F_{\lambda}\right)\right]
\eee
We now integrate by parts to compute using \fref{defpotential}:
\bee
\int\frac{b(\Lambda Z)_\l}{\l^2r}w_3 A_\l^*w_3 & = & \frac{b}{\l^8}\int \frac{\Lambda Z}{y}\e_3A^*\e_3=\frac{b}{\l^8}\int\frac{2(1+Z)\Lambda Z-\Lambda^2Z}{2y^2}\e_3^2\\
& = & \frac{b}{\l^8}\int\frac{\Lambda\tilde{V}}{2y^2}\e_3^2=b \int \frac{(\Lambda \tilde{V})_\l}{2\l^2r^2} w_3^2.
\eee
Injecting this into the energy identity \fref{firstestimate} yields the modified energy identity:
\bea
\label{modifeideenrgy}
\nonumber &&\frac{1}{2}\frac{d}{dt}\left\{ \mathcal E_4+2\int \frac{b(\Lambda Z)_\l}{\l^2r}w_3w_2\right\} =  - \int (\tilde{H}_\l w_3)^2\\
\nonumber & - & \left(\lsl+b\right)\int \frac{(\Lambda \tilde{V})_\l}{2\l^2r^2} w_3^2+\int \frac{d}{dt}\left(\frac{b(\Lambda Z)_\l}{\l^2r}\right)w_3w_2\\
\nonumber & + & \int \tilde{H}_\l w_3\left[\frac{\partial_{t} Z_{\lambda}}{r} w_2 + \int A_{\lambda}\left(\frac{\partial_{t} V_{\lambda}}{r^2}w \right)+\left(AH \right)_{\lambda}\left( \frac{1}{\lambda^2}\mathcal F_{\lambda}\right)\right]\\
\nonumber & + & \int \frac{b(\Lambda Z)_\l}{\l^2r}w_2\left[-\tilde{H}_\l w_3+\frac{\partial_{t} Z_{\lambda}}{r} w_2 + A_{\lambda} \left( \frac{\partial_{t} V_{\lambda}}{r^2} w \right) + \left(AH \right)_{\lambda} \left(\frac 1{\lambda^2} \mathcal F_{\lambda}\right)\right]\\
& + & \int\frac{b(\Lambda Z)_\l}{\l^2r}w_3\left[\frac{\partial_{t} V_{\lambda}}{r^2} w + H_{\lambda} \left(\frac 1{\lambda^2} \mathcal F_{\lambda}\right)\right]
\eea
We now aim at estimating all terms in the RHS of \fref{modifeideenrgy}. All along the proof, we shall make an implicit use of the coercitivity estimates of Lemma \ref{coerchtilde} and Lemma \ref{lemmainterpolation}.\\

{\bf step 3} Lower order quadratic terms.\\

We start with treating the lower order quadratic terms in \fref{modifeideenrgy} {\it using dissipation}. Indeed, we have from \fref{comportementz}, \fref{comportementv}, \fref{rougboundpope}  the bounds:
\bea
|\partial_t Z_{\lambda}| + |\partial_t V_{\lambda}| \lesssim \frac{b}{\lambda^2} \left( |\Lambda Z| + |\Lambda V|\right)_{\lambda} \lesssim \frac{b}{\lambda^2}\frac{y^2}{1+y^4}
\eea
and thus from Cauchy Schwartz, the rough bound \fref{rougboundpope} and Lemma \ref{lemmainterpolation}:
\bee
&&\int\left|\tilde{H}_\l w_3\left[\frac{\partial_{t} Z_{\lambda}}{r} w_2 + \int A_{\lambda}\left(\frac{\partial_{t} V_{\lambda}}{r^2}w \right)\right]\right|+\int|\tilde{H}_\l w_3|\left|\frac{b(\Lambda Z)_\l}{\l^2r}w_2\right|\\
& \leq &\frac{1}{2}\int|\tilde{H}_\l w_3|^2+\frac{b^2}{\l^8}\left[\int \frac{\varepsilon_2^2}{1+y^6}+\int \frac{\varepsilon_1^2}{1+y^8} + \int \frac{\varepsilon^2}{y^2(1+y^8)}\right]\\
& \leq &\frac{1}{2}\int|\tilde{H}_\l w_3|^2+\frac{b}{\l^8}b|\log b|^C b^4.
\eee
All other quadratic terms are lower order by a factor $b$ using again \fref{rougboundpope} and Lemma \ref{lemmainterpolation}:
\bee
&&\left|\lsl+b\right|\int \left|\frac{(\Lambda \tilde{V})_\l}{2\l^2r^2} w_3^2\right|+\int\left|\frac{b(\Lambda Z)_\l}{\l^2r}w_2\left[\frac{\partial_{t} Z_{\lambda}}{r} w_2 + A_{\lambda} \left( \frac{\partial_{t} V_{\lambda}}{r^2} w \right)\right]\right|\\
& + & \int\left|\frac{b(\Lambda Z)_\l}{\l^2r}w_3\frac{\partial_{t} V_{\lambda}}{r^2} w\right|+\left|\int \frac{d}{dt}\left(\frac{b(\Lambda Z)_\l}{\l^2r}\right)w_3w_2\right|\\
& \lesssim & \frac{b^2}{\l^8}\left[\int\frac{\e_3^2}{1+y^2}+\int \frac{\varepsilon_2^2}{1+y^4}+\int \frac{\varepsilon_1^2}{1+y^8} + \int \frac{\varepsilon^2}{y^2(1+y^8)}\right]\\
& \lesssim & \frac{b}{\l^8}b|\log b|^Cb^4.
\eee 
We similarily estimate the boundary term in time: $$\left|\int \frac{b(\Lambda Z)_\l}{\l^2r}w_3w_2\right|\lesssim \frac{b}{\l^6}\left[\int\frac{\e_3^2}{1+y^2}+\int\frac{\e_2^2}{1+y^4}\right]\lesssim \frac{b}{\l^6}|\log b|^Cb^4.$$
We inject these estimates into \fref{modifeideenrgy} to derive the preliminary bound:
\bea
\label{neoheohohe}
 &&\frac{1}{2}\frac{d}{dt}\left\{ \frac{1}{\l^6}\left[\mathcal E_4+O\left(\sqrt{b}\frac{b^4}{|\log b|^2}\right)\right]\right\}\leq -\frac12\int(\tilde{H}_\l w_3)^2+\frac{b}{\l^8}\frac{b^4}{|\log b|^2}\\
\nonumber & + & \int \tilde{H}_\l w_3A_\lambda H_\l\left( \frac{1}{\lambda^2}\mathcal F_{\lambda}\right)+\int H_{\lambda} \left(\frac 1{\lambda^2} \mathcal F_{\lambda}\right)\left[\frac{b(\Lambda Z)_\l}{\l^2r}w_3+A^*_{\l}\left( \frac{b(\Lambda Z)_\l}{\l^2r}w_2\right)\right]
\eea
with constants independent of $M$ for $|b|<b^*(M)$ small enough. We now aim at estimating all terms in the RHS of \fref{neoheohohe}.\\

{\bf step 4} Further use of dissipation.\\

Let us introduce the decomposition from \fref{eqepsilon}, \fref{defF}:
$$\mathcal F=\mathcal F_0+\mathcal F_1, \ \ \mathcal F_0=-\Psit_b-Mod(t), \ \ \mathcal F_1=L(\e)-N(\e).$$
The first term in the RHS of \fref{neoheohohe} is estimated after an integration by parts:
\bea
\label{oneone}
\nonumber&& \left| \int \tilde{H}_\l w_3A_\lambda H_\l\left( \frac{1}{\lambda^2}\mathcal F_{\lambda}\right)\right|\leq \frac{C}{\l^8}\|A^*\e_3\|_{L^2}\|H^2\mathcal F_0\|_{L^2}+\frac14\int(\tilde{H}_\l w_3)^2+\frac{C}{\l^8}\int|AH\mathcal F_1|^2\\
& \leq & \frac{C}{\l^8}\left[\|H^2\mathcal F_0\|_{L^2}\sqrt{\mathcal E_4}+\|AH\mathcal F_1\|_{L^2}^2\right]+\frac14\int(\tilde{H}_\l w_3)^2
\eea
for some universal constant $C>0$ independent of $M$. The last two terms in \fref{neoheohohe} can be estimated in brute force from Cauchy Schwarz:
\bea
\label{onetwo}
\nonumber \left|\int H_{\lambda} \left(\frac 1{\lambda^2} \mathcal F_{\lambda}\right)\frac{b(\Lambda Z)_\l}{\l^2r}w_3\right| & \lesssim&  \frac{b}{\l^8}\left(\int\frac{1+|\log y|^2}{1+y^4}|H\mathcal F|^2\right)^{\frac12}\left(\int\frac{\e_3^2}{y^2(1+|\log y|^2)}\right)^{\frac12}\\
& \lesssim & \frac{b}{\l^8}\sqrt{\mathcal E_4}\left(\int\frac{1+|\log y|^2}{1+y^4}|H\mathcal F|^2\right)^{\frac12}
\eea
where constants are independent of $M$ thanks to the estimate \fref{coercwthree} for $\e_3$. Similarily:
\bea
\label{onethree}
\nonumber &&\left|\int H_{\lambda} \left(\frac 1{\lambda^2} \mathcal F_{\lambda}\right)A^*_{\l}\left( \frac{b(\Lambda Z)_\l}{\l^2r}w_2\right)\right|\\
\nonumber &  \lesssim &  \frac{b}{\l^8}\left(\int\frac{1+|\log y|^2}{1+y^2}|AH\mathcal F|^2\right)^{\frac12}\left(\int\frac{\e_2^2}{1+y^4(1+|\log y|^2)}\right)^{\frac12}\\
 & \lesssim & \frac{b}{\l^8}C(M)\sqrt{\mathcal E_4}\left(\int\frac{1+|\log y|^2}{1+y^2}|AH\mathcal F_0|^2+\int|AH\mathcal F_1|^2\right)^{\frac12}.
\eea
We now claim the bounds:
\be
\label{crucialboundthree}
\int\frac{1+|\log y|^2}{1+y^4}|H\mathcal F|^2\lesssim \frac{b^4}{|\log b|^2}+\frac{\mathcal E_4}{\log M},
\ee
\be
\label{weigheivbiovheo}
\int\frac{1+|\log y|^2}{1+y^2}|AH\mathcal F_0|^2 \lesssim \delta^*\left(\frac{b^4}{|\log b|^2}+\mathcal E_4\right),
\ee
\be
\label{cnofooeeo}
\int|H^2\mathcal F_0|^2\lesssim b^2\left[\frac{b^4}{|\log b|^2}+\frac{\mathcal E_4}{\log M}\right],
\ee
\be
\label{crucialboundtwo}
\int |AH\mathcal F_1|^2\lesssim b\left[\frac{b^4}{|\log b|^2}+\frac{\mathcal E_4}{\log M}\right],
\ee
with all $\lesssim$ constants independent of $M$ for $|b|<b^*(M)$ small enough, and where $$\delta^*=\delta^*(b^*(M))\to 0\ \ \mbox{as} \ \ b^*(M)\to 0.$$ Injecting these bounds together with \fref{oneone}, \fref{onetwo}, \fref{onethree} into \fref{neoheohohe} concludes the proof of \fref{monoenoiencle}. We now turn to the proof of \fref{crucialboundthree}, \fref{weigheivbiovheo}, \fref{cnofooeeo}, \fref{crucialboundtwo}.\\

{\bf step 5} $\Psit_b$ terms.\\

The contribution of $\Psit_b$ terms to \fref{crucialboundthree}, \fref{weigheivbiovheo}, \fref{cnofooeeo} is estimated from \fref{wiehgthtwoewighttilde}, \fref{wiehgthtwoewightthreetilde}, \fref{controleh4erreurtilde} which are at the heart of the construction of $\tilde{Q}_b$ and yield the desired bounds.\\

{\bf step 6} $Mod(t)$ terms.\\

Recall the definition \fref{defmode} of $Mod(t)$:
$$Mod(t)=- \left( \frac{\lambda_s}{\lambda} +b\right) \Lambda \qbt + (b_s + b^2)(\tilde{T_1} +2b\tilde{T_2}).$$ For \fref{crucialboundthree}, we estimate using the rough bounds \fref{roughtwto}, \fref{roughtwtobis} and the control of the modulation parameters \fref{parameters}, \fref{parameterspresicely} to estimate:
\bee
&&\int\frac{1+|\log y|^2}{1+y^4}|H Mod|^2\\
& \lesssim &\left|\lsl+b\right|^2\int\frac{1+|\log y|^2}{1+y^4}|H\Lambda \qbt|^2+|b_s+b^2|^2\int\frac{1+|\log y|^2}{1+y^4}|H(\tt_1+2b\tt_2)|^2\\
&\lesssim & \left|\lsl+b\right|^2+|b_s+b^2|^2\lesssim \frac{b^4}{|\log b|^2}+\frac{\mathcal E_4}{\log M}.
\eee
For \fref{weigheivbiovheo}, we use the cancellations $H\Lambda Q=0$, $AHT_1=0$ and the rough bounds \fref{roughtwto}, \fref{roughtwtobis} to derive the degenerate bounds:
\bee
\int\frac{1+|\log y|^2}{1+y^2}|AH\Lambda \qbt|^2&\lesssim& \int_{y\leq 2B_1}\frac{1+|\log y|^2}{1+y^2}\left|\frac{by|\log y|+b^2(1+y^3)+b^3(1+y^5)}{1+y^3}\right|^2\\
& \lesssim & b^2
\eee
\bee
\int\frac{1+|\log y|^2}{1+y^2}|AH(\tt_1+b\tt_2)|^2& \lesssim& \int\frac{1+|\log y|^2}{1+y^2}\left|\frac{y\log y}{y^3}{\bf 1}_{B_1\leq y\leq 2B_1}+b\frac{1+y^3}{1+y^3}{\bf 1}_{y\leq 2B_1}\right|^2\\
& \lesssim & b
\eee
and thus from \fref{parameters}, \fref{parameterspresicely}:
\bee 
\int\frac{1+|\log y|^2}{1+y^2}|AH\mathcal Mod|^2\lesssim b\left[\left(\lsl+b\right)^2+(b_s+b^2)^2\right]\lesssim b\left[\frac{b^4}{|\log b|^2}+\frac{\mathcal E_4}{\log M}\right].
\eee
For \fref{cnofooeeo}, we estimate from the rough bounds \fref{roughtwto}, \fref{roughtwtobis}:
$$
\int|H^2\Lambda \qbt|^2\lesssim \int \left|\frac{by\log y+b^2(1+y^3)+b^3(1+y^5)}{1+y^4}\right|^2\lesssim b^2,
$$
$$\int|H^2\tt_1|^2\lesssim \int_{B_1\leq y\leq 2B_1}\left|\frac{y\log y}{y^4}\right|^2\lesssim \frac{|\log b|^2}{B_1^4}\lesssim b^2.$$
The last term is more subtle and we claim:
\be
\label{remnaininini}
\int|H^2\tt_2|^2\lesssim 1
\ee
which yields
$$\int|H^2Mod|^2\lesssim b^2\left[\left(\lsl+b\right)^2+(b_s+b^2)^2\right]\lesssim b^2\left[\frac{b^4}{|\log b|^2}+\frac{\mathcal E_4}{|\log M}\right]$$ and concludes the proof of \fref{cnofooeeo}.\\
{\it Proof of \fref{remnaininini}}: First by definition of $\tt_2=\chi_{B_1}T_2$, the rough bound \fref{roughtwto} and \fref{eqtzerotwo}:
\be
\label{cpeopeepi}
\int|H^2\tt_2|^2\lesssim \left[\int_{B_1\leq y\leq 2B_1}\left|\frac{y^3}{y^4}\right|^2+\int_{y\leq 2B_1}|H\Sigma_2|^2\right]\lesssim 1+\int_{y\leq 2B_1}|H\Sigma_2|^2.
\ee We now compute from \fref{defsigma2}, \fref{defsigmab}:
\bee
H\Sigma_2 & = &  H\Sigma_b + H(T_1 - \Lambda T_1) +O\left(\frac{y^2|\log y|^2}{1+y^5}\right)\\
& = & -c_b\chi_{\frac{B_0}{4}}\Lambda Q +d_b H[(1 - \chi_{3B_0}) \Lambda Q]+ H(T_1 - \Lambda T_1) +O\left(\frac{y^2|\log y|^2}{1+y^5}\right)\\
& = & \frac{1}{|\log b|}O\left(\frac{1}{1+y}{\bf1}_{y\leq 3B_0}\right)+H(T_1 - \Lambda T_1) +O\left(\frac{y^2|\log y|^2}{1+y^5}\right)
\eee
and observe using \fref{calculrealcinr}, the asympotics \fref{infinity} of $\Lambda Q$ and the fundamental cancellation $\left(\Lambda+\Lambda^2\right)\left(\frac1y\right)=0$ that:
\bee
HT_1-H\Lambda T_1 & = & HT_1-\left(2HT_1+\Lambda HT_1-\frac{\Lambda V}{y^2}T_1\right)=\Lambda Q+\Lambda^2Q+O\left(\frac{y\log y}{1+y^4}\right)\\
& = & O\left(\frac{\log y}{1+y^3}\right).
\eee
We thus conclude: $$\int |H\Sigma_2|^2\lesssim \frac{1}{|\log b|^2}\int_{y\leq 2B_0}\frac{1}{1+y^2}+\int_{y\leq 2B_1}\frac{|\log y|^4}{1+y^6}\lesssim 1$$ which together with \fref{cpeopeepi} concludes the proof of \fref{remnaininini}.\\

{\bf step 7} Small linear term $L(\e)$.\\

Let us rewrite from  a Taylor expansion: 
\be
\label{cnoheiohoe}
L(\e)=-\frac{N_2(\tilde{\alpha})}{y^2}\e\ \ \mbox{with} \ \ N_2(\tilde{\alpha})=f'(Q+\tilde{\alpha})-f'(Q)=\tilde{\alpha}\int_0^1f''(Q+\tau \alphat)d\tau.
\ee
Near the origin $y\leq 1$, we use $f''(0)=0$ and the estimate $|\alphat|\lesssim by^3$ near the origin by construction to obtain the high order cancellation 
\be
\label{estimateorigiinr}
N_2(\tilde{\alpha})\lesssim |\alphat|y\lesssim by^4
\ee 
which together with the bounds \fref{estoriiginagian}, \fref{choeouefeouie}, \fref{estlinftydeux} easily yields: $$\int_{y\leq 1}\frac{1+|\log y|^2}{1+y^4}|H L(\e)|^2+\int_{y\leq 1}|AHL(\e)|^2\lesssim b^6.$$ For $y\geq 1$, we use $f''(\pi)=0$ and the bounds \fref{decayttwo}, \fref{T3} to derive the bound: 
\bea
\label{boundjjs}
\nonumber |N_2(\alphat)| & \lesssim & |\alphat|\left[\frac{1}{1+y}+|\alphat|\right]\lesssim \left[b|\log y|+b^2y^2|\log y|^2\right]{\bf 1}_{y\leq 2B_1}\\
& \lesssim & b|\log b|^C{\bf 1}_{y\leq 2B_1}.
\eea 
A brute force computation taking futher derivatives and using \fref{estun} now yields the control: 
\bee
\int_{y\geq 1}|AHL(\e)|^2+\int_{y\geq 1}\frac{1+|\log y|^2}{1+y^4}|H L(\e)|^2 &\lesssim &b^2|\log b|^C\int_{y\leq 2B_1}\Sigma_{i=0}^3\frac{|\pa^i_y\e|^2}{y^{2(4-i)}}\\
&\lesssim &b^2|\log b|^C\mathcal E_4\lesssim  b\delta^*\mathcal E_4.
\eee
This concludes the proof of \fref{crucialboundthree}, \fref{crucialboundtwo} for $L(\e)$ terms.\\

{\bf step 8} Nonlinear term $N(\e)$.\\

Let us now treat the nonlinear term \fref{defNe}. We split the contribution at the origin and far out and claim:
\be
\label{neattheorigin}
\forall y\leq 1, \ \  |AHN(\e)(y)|\lesssim C(M)b^4|\log y|^4,
\ee
\be
\label{estimatefarout}
\int_{y\geq 1} (AHN(\e))^2\lesssim C(M)\frac{b^6}{|\log b|^2}.
\ee
which implies \fref{crucialboundtwo} for $N(\e)$. The estimate \fref{crucialboundthree} follows along similar lines and is in fact simpler and left to the reader.\\
{\it Proof of  \fref{neattheorigin}}: We need to treat the possible singularity at the origin. For this, let us rewrite by Taylor expansion: $$N(\e)=z^2N_0(\e)\ \ \mbox{with} \ \ z=\frac{\e}{y}, \ \ N_0(\e)=\int_0^1(1-\tau)f''(\qbt+\tau \e)d\tau,$$ and thus: $$H(N(\e))=-N_0(\e)\Delta(z^2)-2z\pa_yz\pa_yN_0(\e)+z^2H(N_0(\e)),$$ 
\bea
\label{formaula}
AHN(\e)& = & -A(N_0(\e))\Delta(z^2)+N_0(\e)\pa_y\Delta(z^2)-2A(z\pa_yz)\pa_y(N_0(\e))\\
\nonumber & + & 2z\pa_yz\pa_{yy}(N_0(\e))+z^2AH (N_0(\e))-2z\pa_yzH(N_0(\e)).
\eea
We now use the Hardy bounds \fref{estoriiginagian}, \fref{choeouefeouie}, \fref{estlinftydeux} and the degeneracy $|Z-1|+|V-1|\lesssim y^2$ for $y\leq 1$ to estimate for $0<y\le 1$:
$$|\pa_{yy}\e|=\left|-H\e+\frac{A\e}{y}+\frac{1-Z+V-1}{y^2}\e\right|\lesssim b^2y|\log y|,$$
$$|z|=\left|\frac{\e}{y}\right|\lesssim b^2,$$
$$|\pa_yz|=\left|-\frac{A\e}{y}+\frac{Z-1}{y^2}\e\right|\lesssim b^2y|\log y|,$$
$$|\pa_{yy}z|=\left|\frac{2+Z}{y^2}A\e-\frac{H\e}{y}+\pa_y\left(\frac{Z-1}{y^2}\e\right)\right|\lesssim b^2|\log y|,$$
\bee
|\pa_{yyy}z| & = & \left|\pa_y\left(\frac{2+Z}{y^2}\right)A\e+\frac{2+Z}{y^2}\left(H\e-\frac{1+Z}{y}A\e\right)+\frac{H\e}{y^2}\right.\\
& - & \left .\frac{1}{y}\left(-AH\e+\frac{Z}{y}H\e\right)+\pa_{yy}\left(\frac{Z-1}{y^2}\e\right)\right|\\
& \lesssim & \frac{b^2|\log y|}{y}
\eee
with constants depending on $M$. We now estimate using $f^{(2k)}(0)=0$: for $0<y\leq 1$,
$$|N_0(\e)|\lesssim y,$$ $$|\pa_yN_0(\e)|=\left|\int_0^{1}(1-\tau)\pa_y(\qbt+\tau\e)f^{(3)}(\qbt+\tau\e)d\tau\right|\lesssim 1,$$  
\bee
|\pa_{yy}N_0(\e)| & = & \left|\int_0^{1}(1-\tau)\left[\pa_{yy}(\qbt+\tau\e)f^{(3)}(\qbt+\tau\e)+\left(\pa_y(\qbt+\tau\e)\right)^2f^{(4)}(\qbt+\tau\e)\right]d\tau\right|\\
& \lesssim & y|\log y|^2
\eee
We need to exploit further cancellations for $AN_0(\e)$ and this requires pushing the Taylor expansion: 
$$N_0(\e)=\frac{1}{2}f''(\qbt)+\e N_1(\e), \ \ N_1(\e)=\int_0^1\int_0^1(1-\sigma)\tau(1-\tau)f^{(3)}(\qbt+\sigma\tau\e)d\sigma d\tau.$$ By construction, $\qbt$ is a smooth function at the origin and admits a Taylor expansion $$\qbt=c_1(b)y+c_2(b)y^3+O(y^5)\ \ \mbox{with} \ \ |c_1(b), c_2(b)|\lesssim 1$$ and hence:
$$|Af''(\qbt)|+|Hf''(\qbt)|+|AH(f''(\qbt))|\lesssim 1.$$ We therefore estimate arguing like for $N_0(\e)$: $$|A(\e N_1(\e))|=|A\e N_1(\e)-\e\pa_y(N_1(\e))|\lesssim y^2|\log y|,$$ $$|H(\e N_1(\e))|=|N_1(\e)H\e-2\pa_y\e\pa_yN_1(\e)-\e\pa_{yy}N_1(\e)|\lesssim y|\log y|^2,$$ 

\bee
|AH(\e N_1(\e))|& = & \left|N_1(\e)AH\e-(H\e)\pa_yN_1(\e)-2(A\pa_y\e)\pa_yN_1(\e)+2\pa_y\e\pa_{yy}N_1(\e)\right.\\
& - & \left.A\e\pa_{yy}N_1(\e) +\e\pa_{yyy}N_1(\e)\right|\\
& \lesssim & |\log y|^4
\eee
Injecting the collection of above estimates into \fref{formaula} now yields \fref{neattheorigin}.\\
{\it Proof of \fref{estimatefarout}}: For $y\geq 1$, we estimate from \fref{estlinftydeuxfa}, \fref{estlinftysecond}: $$\|z\|_{L^{\infty}(y\geq 1)}\lesssim b|\log b|^C, \ \ \|\pa_yz\|_{L^{\infty}(y\geq 1)}+\|\frac{z}{y}\|_{L^{\infty}(y\geq 1)}\lesssim b^{\frac32}|\log b|^C.$$ The construction of $\qbt$ yields the bounds for $y\geq 1$:
$$|\pa_y\qbt|\lesssim |\log b|^C\left(\frac{1}{y^2}+b{\bf 1}_{y\leq 2B_1}\right),$$ $$|\pa_{yy}\qbt|\lesssim|\log b|^C\left( \frac{1}{y^3}+\frac{b}{y}{\bf 1}_{y\leq 2B_1}\right)\lesssim |\log b|^C\left(\frac{1}{y^3}+b^{\frac32}\right),$$ $$|\pa_{yyy}\qbt|\lesssim |\log b|^C\left(\frac{1}{y^4}+\frac{b}{y^2}{\bf 1}_{y\leq 2B_1}\right)\lesssim |\log b|^C\left(\frac{1}{y^4}+b^{2}\right)$$
 which together with \fref{estlinftydeuxfa}, \fref{estlinftysecond}, \fref{estoeofj} yields the pointwise bounds:
 \be
 \label{boundnoone}
 |N_0(\e)|\lesssim 1,
 \ee
 \be
 \label{boundnotwo}
  |\pa_yN_0(\e)|\lesssim |\log b|^C\left(\frac{1}{y^2}+b+\|\pa_y\e\|_{L^{\infty}(y\geq 1)}\right)\lesssim |\log b|^C\left(\frac{1}{y^2}+b\right),
  \ee
 \be
 \label{boundnothree}
 |\pa_{yy}N_0(\e)|\lesssim |\log b|^C\left[\frac{1}{y^3}+b^{\frac32}+\|\pa_{yy}\e\|_{L^{\infty}(y\ge 1)}\right]\lesssim  |\log b|^C\left[\frac{1}{y^3}+b^{\frac32}\right],
 \ee
 \be
 \label{boundnofour}
 |\pa_{yyy}N_0(\e)|\lesssim |\log b|^C\left[\frac{1}{y^4}+b^2+\|\pa_{yyy}\e\|_{L^{\infty}(y\geq 1)}\right]\lesssim |\log b|^C\left[\frac{1}{y^4}+b^2\right].
 \ee
 We now compute:
\be
\label{formulahne}
H(N(\e))=N_0(\e)H(z^2)-2z\pa_yz\pa_yN_0(\e)-z^2\Delta(N_0(\e)),
\ee 
\bee
AHN(\e)& = & N_0(\e)AH(z^2)-\pa_yN_0(\e)H(z^2)-2A(z\pa_yz)\pa_y(N_0(\e))\\
\nonumber & + & 2z\pa_yz\pa_{yy}(N_0(\e))-A(z^2)\Delta(N_0(\e))+z^2\pa_y\Delta (N_0(\e)),
\eee
and hence using the $L^2$ weighted bounds \fref{estun}, \fref{lossyboundwperp}, \fref{inteproloatedbound} and the $L^{\infty}$ bounds \fref{estlinftydeuxfa}, \fref{estlinftysecond}:
\bee
\int_{y\ge 1} |N_0(\e)AH(z^2)|^2 & \lesssim & \int_{y\ge 1}\left[|\pa_yz\Delta z|^2+|z\pa_y\Delta z|^2+|\pa_yz\pa_{yy}z|^2+\frac{|z\pa_yz|^2}{y^4}\right.\\
&+ & \left. \frac{1}{y^2}\left(|z\Delta z|^2+|\pa_yz|^4+\frac{|z|^4}{y^4}\right)\right]\\
& \lesssim & b^6|\log b|^C,
\eee
\bee
\int_{y\ge 1}|\pa_yN_0(\e)H(z^2)|^2\lesssim \int_{y\geq 1} \left(\frac{1}{y^2}+b^2|\log b|^C\right)\left(|z|^2|\Delta z|^2+|\pa_yz|^4+\frac{|z|^4}{y^4}\right)\lesssim b^6|\log b|^C,
\eee
\bee
\int_{y\ge 1}|A(z\pa_yz)\pa_yN_0(\e)|^2& \lesssim &\int_{y\geq 1} \left(\frac{1}{y^2}+b^2|\log b|^C\right)\left(|z|^2|\pa_{yy}z|^2+|\pa_yz|^4+\frac{|z|^2|\pa_yz|^2}{y^2}\right)\\
& \lesssim & b^6|\log b|^C,
\eee
$$
\int_{y\ge 1} |z\pa_yz\pa_{yy}N_0(\e)|^2\lesssim |\log b|^C \int_{y\ge 1}|z|^2|\pa_yz|^2\left[\frac{1}{y^6}+b^3\right]\lesssim b^6|\log b|^C,
$$
$$\int_{y\ge 1}|A(z^2)\Delta N_0(\e)|^2\lesssim |\log b|^C\int_{y\ge 1}\left[|z|^2|\pa_yz|^2+\frac{|z|^4}{y^2}\right]\left[\frac{1}{y^6}+b^3\right]\lesssim b^6|\log b|^C,$$
$$\int_{y\ge 1} |z^2\pa_y\Delta (N_0(\e))|^2\lesssim |\log b|^C\int_{y\ge 1} |z|^4\left[\frac{1}{y^8}+b^4\right]\lesssim b^6|\log b|^C.$$
This concludes the proof of \fref{estimatefarout}.\\
This concludes the proof of \fref{crucialboundthree}, \fref{weigheivbiovheo}, \fref{cnofooeeo}, \fref{crucialboundtwo} and thus of Proposition \ref{AEI2}.

  
  \section{Sharp description of the singularity formation}
  

In this section, we start with completing the proof of the bootstrap Proposition \ref{bootstrap}. Theorem \ref{thmmain} will then easily follow.


\subsection{Closing the bootstrap}
\label{sectionbootstrap}

We are now in position to close the boostrap bounds of Proposition \ref{bootstrap}.\\

{\bf Proof of Proposition \ref{bootstrap}}\\

{\bf step 1} Energy bound.\\

First observe that \fref{init3hbis} and the modulation equation \fref{parameterspresicely} ensure for $M$ large enough\footnote{recall that $K$ in the bootstrap bounds is large but independent of $M$.}: $$b_s\leq 0$$ and the upper bound in \fref{init1hb} follows. We now claim: 
$$b(t)>0\ \ \mbox{on}\  \ [0,T_1).$$ Indeed, from \fref{rougboundpope}, if $b(t_0)=0$ for some $t_0\in [0,T_1)$, then $b(t)\equiv 0$ on some $[t_0-\delta,t_0]$ and thus from \fref{rougboundpope}, \fref{init3hbis}, $\lambda(t)\equiv \lambda(t_0)$ and $u(t)\equiv Q_{\lambda(t_0)}$ on $[t_0-\delta,t_0]$. We can thus iterate on $\delta>0$ and conclude that $u(0)$ is initially a harmonic map, a contradiction. This concludes the proof of \fref{init1hb}.\\
We now prove \fref{init3hb} which follows from the conservation of energy. Indeed, let $$\tilde{\e}=\e+\alphat,$$ then
\bea
\label{consebfie}
E_0 & = & \int|\pa_y(Q+\et)|^2+\int\frac{g^2(Q+\et)}{y^2}\\
\nonumber & = & E(Q)+(H\e,\e)+\int\frac{1}{y^2}\left[g^2(Q+\et)-2f(Q)\e-f'(Q)\et^2\right].
\eea
We now recall from Lemma \ref{coerchtilde} and Lemma \ref{lemmacoer} the coercivity properties: $$(H\e,\e)\geq c(M)\left[\int|\pa_y\e|^2+\int\frac{|\e|^2}{y^2}\right],$$
$$(H\et,\et)\geq c(M)\left[\int|\pa_y\et|^2+\int\frac{|\et|^2}{y^2}\right]+O(C(M)b^2).$$
The nonlinear term is estimated from a Taylor expansion:
$$
\left|\int\frac{1}{y^2}\left[g^2(Q+\et)-2f(Q)\et-f'(Q)\et^2\right]\right|\lesssim  \int \frac{|\et|^3}{y^2}\lesssim  \left(\int|\pa_y\et|^2+\int\frac{|\et|^2}{y^2}\right)^{\frac{3}{2}}.
$$
where we used the Sobolev bound $$\|\et\|_{L^{\infty}}^2\lesssim \|\pa_y\et\|_{L^2}\|\frac{\et}{y}\|_{L^2}.$$
We inject these bounds into the conservation of energy \fref{consebfie} and use the bound on the profile $$\int|\pa_y\alphat|^2+\int\frac{|\alphat|^2}{y^2}\lesssim b|\log b|^C$$ and \fref{init2}, \fref{init1hb} to estimate: 
\bee
\int|\pa_y\e|^2+\int\frac{|\e|^2}{y^2} & \lesssim & \int|\pa_y\et|^2+\int\frac{|\et|^2}{y^2}+b|\log b|^C\lesssim C(M)\left|E_0-E(Q)\right|+O(C(M)b)\\
& \lesssim & C(M)b(0)\leq \sqrt{b(0)}
\eee
for $|b(0)|\leq b^*(M)$ small enough, and \fref{init3hb} is proved.\\

{\bf step 2} Control of $\mathcal E_4$.\\ 
 
 We now close the bootstrap bound \fref{init3hbisb} which follows by reintegrating the Lyapounov monotonicity \fref{monoenoiencle} in the regime governed by the modulation equations \fref{parameters}, \fref{parameterspresicely}. Indeed, inject the bootstrap bound \fref{init3hbis} into the monotonicity formula \fref{monoenoiencle} and integrate in time; this yields: $\forall t\in [0,T_1)$,
\bea
\label{monotnyintegree}
\mathcal E_4(t) & \leq & 2\left(\frac{\lambda(t)}{\lambda(0)}\right)^6\left[\mathcal E_4(0)+C\sqrt{b(0)}\frac{b^4(0)}{|\log b(0)|^2}\right]+\frac{b^4(t)}{|\log b(t)|^2}\\
\nonumber & + & C\left[1+\frac{K}{\log M}+\sqrt{K}\right]\lambda^6(t)\int_0^t\frac{b}{\lambda^8}\frac{b^4}{|\log b|^2}d\tau
\eea
for some universal constant $C>0$ independent of $M$.\\
Let us now consider two constants 
\be
\label{conno}
\alpha_1=2-\frac{C_1}{\sqrt{\log M}}, \ \ \alpha_2=2+\frac{C_2}{\sqrt{\log M}}
\ee for some large enough universal constanst $C_1,C_2$. We compute using the modulation equations \fref{parameters}, \fref{parameterspresicely} and the bootstrap bound \fref{init3hbis}:
\bee
\frac{d}{ds}\left\{\frac{|\log b|^{\alpha_i}b}{\lambda}\right\}& = & \frac{|\log b|^{\alpha_i}}{\lambda}\left[\left(1-\frac{\alpha_i}{|\log b|}\right)b_s-\frac{\lambda_s}{\lambda}b\right]\\
& = & \frac{|\log b|^{\alpha_i}}{\lambda}\left[\left(1-\frac{\alpha_i}{|\log b|}\right)b_s+b^2 + O\left(\dfrac{b^3}{|\log b|}\right)\right]\\
& = & \left(1-\frac{\alpha_i}{|\log b|}\right)\frac{|\log b|^{\alpha_i}}{\lambda}\left[b_s+b^2\left(1+\frac{\alpha_i}{|\log b|}+O\left(\frac{1}{|\log b|^2}\right)\right)\right]\\
& &\left\{ \begin{array}{ll}\leq 0 \ \ \mbox{for} \ \ i=1\\\geq 0 \ \ \mbox{for}  \  \ i=2.
\end{array}\right .
\eee
Integrating this from $0$ to $t$ yields: 
\be
\label{lawintegrationone}
\frac{b(0)}{\lambda(0)}\left(\frac{|\log b(0||}{|\log b(t)|}\right)^{\alpha_2}
\leq \frac{b(t)}{\lambda(t)}\leq \frac{b(0)}{\lambda(0)}\left(\frac{|\log b(0||}{|\log b(t)|}\right)^{\alpha_1}.
\ee
This yields in particular using the initial bound \fref{init2} and the bound \fref{init1hb}:
\be
\label{estationtermezero}
\left(\frac{\lambda(t)}{\lambda(0)}\right)^6\mathcal E_4(0)\leq (b(t)|\log b(t)|^{\alpha_2})^6\frac{\mathcal E_0}{(b(0)|\log b(0)|^{\alpha_2})^6}\leq \frac{b^4(t)}{|\log b(t)|^2},
\ee 
\bea
\label{estationtermezerobis}
\nonumber C\left(\frac{\lambda(t)}{\lambda(0)}\right)^6\sqrt{b(0)}\frac{b^4(0)}{|\log b(0)|^2} & \lesssim &\left(\frac{b(t)|\log b(t)|^{\alpha_2}}{b(0)|\log b(0)|^{\alpha_2}}\right)^6\sqrt{b(0)}\frac{b^4(0)}{|\log b(0)|^2}\\
& \lesssim & C(b(t))^{4+\frac14}\leq \frac{b^4(t)}{|\log b(t)|^2}.
\eea
We now compute explicitely using $b=-\lambda\lambda_{t}+O\left(\frac{b^2}{|\log b|}\right)$ from \fref{parameters}:
\bee
\int_0^t\frac{b}{\lambda^8}\frac{b^4}{|\log b|^2}d\sigma& = & \frac16\left[\frac{b^4}{\lambda^6|\log b|^2}\right]_0^t-\frac 16\int_0^{t}\frac{b_tb^3}{\lambda^6|\log b|^2}\left(4+\frac{2}{|\log b|}\right)d\tau\\
& + & O\left(\int_0^t\frac{b}{\lambda^8}\frac{b^5}{|\log b|^2}d\tau\right)
\eee
which implies using now $|b_s+b^2|\lesssim \frac{b^2}{|\log b|^2}$ from \fref{parameterspresicely} and \fref{init3hbis}:
$$\lambda^6(t)\int_0^t\frac{b}{\lambda^8}\frac{b^4}{|\log b|^2}d\sigma \lesssim \left[1+O\left(\frac{1}{|\log b_0|}\right)\right]\frac{b^4(t)}{|\log b(t)|^2}.$$ Injecting this together with \fref{estationtermezero}, \fref{estationtermezerobis} into \fref{monotnyintegree} yields $$\mathcal E_4(t)\leq C \frac{b^4(t)}{|\log b(t)|^2}\left[1+\frac{K}{\log M}+\sqrt{K}\right]$$ for some universal constant $C>0$ independent of $K$ and $M$, and thus \fref{init3hbisb}  follows for $K$ large enough independent of $M$.\\

{\bf step 3} Control of $\mathcal E_2$.\\

We now close the $H^2$ bound \fref{init3hb}. This bound is used mostly in the proof of the interpolation estimates of Lemma \ref{lemmainterpolation}, and there the power the log is \fref{init3hb} is irrelevant. It becomes on the contrary critical in the proof of the regularity \fref{regularityustar} and this requires being careful with logarithmic growth. For this reason, the profile $\qbh$ localized near $B_0$ given by \fref{localisation} is better adapted to the $\mathcal E_2$ control.\\
Let then the radiation $\zeta_b$ given by \fref{defraidiation} and the new decomposition of the flow:
\be
\label{dfoeeo}
u=(\qbt+\e)_{\lambda}=(\qbh+\eh)_{\lambda}\ \ \mbox{ie}\ \ \eh=\e+\zeta_b
\ee and the renormalization $$\wh(t,r)=\eh(s,y).$$ The equation for $\wh$ is similarily like \fref{eqenwini}: 
\be
\label{eqwhat}
\pa_t\wh+H_{\l}\wh=\frac{1}{\l^2}\widehat{\mathcal F}_\l, \ \ \widehat{\mathcal F}=-\hat{\Psi}_b-\widehat{Mod}+\hat{L}(\eh)-\hat{N}(\eh),
\ee
$$\hat{L}(\eh)=  \frac{f'(Q)-f'(\qbh)}{y^2}\eh, \ \ \hat{N}(\eh) = \frac{f(\qbh+ \eh) - f(\qbh)-\eh f'(\qbh)}{y^2}.$$
We then let $$\eh_2=H\eh, \ \ \wh_2=H_\lambda \wh,$$ which satisfies from \fref{eqwhat}: 
$$\pa_t\wh_2+H_\l\wh_2=\frac{\pa_tV_\l}{r^2}\wh+H_\l\left(\frac{1}{\l^2}\widehat{\mathcal F}_\l\right).$$ We therefore compute the energy identity:
\bea
\label{cnokcneoeoneo}
 \frac{1}{2}\frac{d}{dt}\left\{\int|\wh_2|^2\right\}& = & \int \wh_2\left[ -H_\l \wh_2+\frac{\partial_{t} V_{\lambda}}{r^2} \wh + H_{\lambda} \left(\frac 1{\lambda^2} \widehat{\mathcal F}_{\lambda}\right)\right]\\
\nonumber& \lesssim & -\int|A_\l \wh_2|^2+\frac{1}{\l^4}\left[b\|\eh_2\|_{L^2}\left\|\frac{\Lambda V}{y^2}\eh\right\|_{L^2}+|(\eh_2,H\widehat{\mathcal F})|\right]
\eea
and aim at estimating all terms in the above RHS. The local term is estimated using the decomposition \fref{dfoeeo}, the estimate \fref{radiation} and Lemma \ref{lemmainterpolation}:
\be
\label{cnkoneoheofeo}
\int\left|\frac{\Lambda V}{y^2}\eh\right|^2\lesssim b^4|\log b|^C+\int\frac{|\zeta_b|^2}{1+y^8}\lesssim b^4|\log b|^C.
\ee
We now claim the bound:
\be
\label{tobepreovehtwo}
|(\eh_2,H\widehat{\mathcal F})|\lesssim b^3|\log b|^2.
\ee
Asssume \fref{tobepreovehtwo}, we then obtain from \fref{cnokcneoeoneo}, \fref{cnkoneoheofeo} the pointwise bound:
$$\frac{d}{dt}\left\{\int|\wh_2|^2\right\}\lesssim \frac{b^3|\log b|^2}{\l^4},$$ which we integrate using also \fref{radiation}:
\bea
\label{cnikheoeif}
\mathcal E_2(t) & = &  \lambda^2(t)\|w_2(t)\|^2_{L^2}\lesssim  \|H\zeta_b(t)\|^2_{L^2}+\lambda^2(t)\|\wh_2(t)\|^2_{L^2}\\
\nonumber & \lesssim & b^4(t)|\log b(t)|^2+\left(\frac{\lambda(t)}{\lambda(0)}\right)^2\left[\mathcal E_2(0)+b^2(0)|\log b(0)|^2\right]+\lambda^2(t)\int_0^t\frac{b^3|\log b|^2}{\lambda^4(\tau)}d\tau.
\eea
From \fref{init2}, \fref{lawintegrationone}:
\bee
\left(\frac{\lambda(t)}{\lambda(0)}\right)^2\left[\mathcal E_2(0)+b^2(0)|\log b(0)|^2\right] & \lesssim & \frac{(b(0))^{10}+b^2(0)|\log b(0)|^2}{(b(0)|\log b(0)|^{\alpha_2}|)^2} b^2(t)|\log b(t)|^{2\alpha_2}\\
& \leq & b^2(t)|\log b(t)|^{4+\frac14},
\eee
We now use the bound $b_s\lesssim -b^2$ and \fref{lawintegrationone} to estimate: 
\bee
&&\lambda^2(t)\int_0^t\frac{b^3|\log b|^2}{\lambda^4(\tau)}d\tau\lesssim  \lambda^2(t)\int_0^t\frac{-b_tb|\log b|^2}{\lambda^2(\tau)}d\tau\\
& \lesssim & \left(\frac{\lambda(t)}{\lambda(0)}\right)^2b^2(0)|\log b(0)|^{2\alpha_1}\int_0^t\frac{-b_t}{b|\log b|^{2\alpha_1-2}}d\tau\\
& \lesssim &  \left(\frac{\lambda(t)}{\lambda(0)}\right)^2b^2(0)|\log b(0)|^{2\alpha_1}\frac{1}{|\log b(0)|^{2\alpha_1-3}}\\
& \lesssim & b^2(t)|\log b(t)|^{2\alpha_2}\frac{|\log b(0)|^3}{|\log b(0)|^{2\alpha_2}}\lesssim b^2(t)|\log b(t)|^{4+\frac 14}.
\eee
Injecting these bounds into \fref{cnikheoeif} yields: $$\mathcal E_2(t)\lesssim b^2(t)|\log b(t)|^{4+\frac14}$$ and concludes the proof of \fref{init3hb}.\\

{\it Proof of \fref{tobepreovehtwo}}: We estimate the contribution of each term in \fref{tobepreovehtwo} coming from the decomposition \fref{eqwhat}. First observe from the interpolation bound \fref{estun} and \fref{radiation}:
\bea
\label{improvedbound}
\nonumber \int_{y\leq 2B_0}|\eh_2|^2 & \lesssim & B_0^4|\log b|^2\int \frac{|\e|^2}{(1+y^4)|\log y|^2}+\int|H\zeta_b|^2\\
& \lesssim & C(M)b^2+b^2|\log b|^2\lesssim b^2|\log b|^2.
\eea
The $\Psih_b$ term is now estimated using  \fref{controleh2erreurtildehat} and \fref{improvedbound}:
\bee
|(\eh_2,H\Psih_b)|& \lesssim & \|H\Psih_b\|_{L^2}\|\eh_2\|_{L^2(y\leq 2B_0)}\lesssim  \left(b^4|\log b|^2b^2|\log b|^2\right)^{\frac12}\lesssim b^3|\log b|^2.
\eee
 We next estimate from \fref{asympT1}, \fref{decayttwo}:
$$ \int|H\hat{T}_1|^2\lesssim \int_{y\leq 2B_0}|\Lambda Q|^2+\int_{B_0\leq y\leq 2B_0}\left|\frac{\log y}{y}\right|^2\lesssim |\log b|^2,$$
$$\int |H\hat{T}_2|^2\lesssim  \int_{y\leq 2B_0}\left|\frac{y}{y^2b|\log b|}\right|^2\lesssim \frac{1}{b^2|\log b|},$$ and thus from \fref{defmodebis}, \fref{parameters}, \fref{parameterspresicely}: 
\bee
\int|H\widehat{Mod}(t)|^2 & \lesssim & \left|\lsl+b\right|^2\int|H\Lambda \qbh|^2+|b_s+b^2|^2\int|H(\hat{T}_1+b\hat{T}_2)|^2\\
& \lesssim & \frac{b^4}{|\log b|^2}|\log b|^2\lesssim b^4|\log b|^2.
\eee
Moreover, $\mbox{Supp}(H\widehat{Mod})\subset[0,2B_0]$ and thus with \fref{improvedbound}:
$$|(\eh_2,H\widehat{Mod})| \lesssim  \left(b^4|\log b|^2b^2|\log b|^2\right)^{\frac12}\lesssim b^3|\log b|^2.
$$
We now turn to the control of the small linear term $\hat{L}(\eh)$ which we rewrite as for \fref{cnoheiohoe}:
$$\hat{L}(\eh)=-\frac{N_2(\alphah)}{y^2}\eh\ \ \mbox{with} \ \ N_2(\alphah)=f'(Q+\alphah)-f'(Q)=\alphah\int_0^1f''(Q+\tau \alphah)d\tau.
$$
Near the origin $y\leq 1$, $\eh=\e$ and the high order vanishing \fref{estimateorigiinr} and the bounds \fref{estoriiginagian}, \fref{choeouefeouie}, \fref{estlinftydeux} easily yield: $$\int_{y\leq 1}|H\hat{L}(\eh)|^2\lesssim b^6.$$ For $y\geq 1$, we estimate like for \fref{boundjjs}: $$ |N_2(\alphah)| \lesssim b|\log b|^C{\bf 1}_{y\leq 2B_0},$$  and then a brute force computation and \fref{radiation} yield the control:
\bee
\int_{y\geq 1}|H\hat{L}(\eh)|^2\lesssim b^2|\log b|^C\int_{y\leq 2B_0}\Sigma_{i=0}^2\frac{|\pa^i_y\eh|^2}{y^{2(4-i)}}\lesssim b^2|\log b|^Cb^4|\log b|^C\lesssim  b^5.
\eee
We also estimate from \fref{radiation} and the bootstrap bound \fref{init3h}: 
\be
\label{htwoehat}
\|\eh_2\|_{L^2}^2\lesssim b^2|\log b|^C
\ee 
and thus $$|(\eh_2,H\hat{L}(\eh))|\lesssim \left(b^2|\log b|^Cb^5\right)^{\frac12}\leq b^3|\log b|^2.$$
It remains to estimate the nonlinear term. Near the origin, we argue like for the proof of \fref{neattheorigin} to derive: $$\forall y\leq 1, \ \ |H(\hat{N}(\eh))|=|H(N(\e))|\lesssim b^4|\log b|^C,$$ this is left to the reader. For $y\geq 1$, we introduce the decomposition:
$$\hat{N}(\eh)=\zh^2\hat{N}_0(\eh)\ \ \mbox{with} \ \ \zh=\frac{\eh}{y}, \ \ \hat{N}_0(\eh)=\int_0^1(1-\tau)f''(\qbh+\tau \eh)d\tau,$$
and recall the formula \fref{formulahne}:
$$H(\hat{N}(\e))=\hat{N}_0(\eh)H(\zh^2)-2\zh\pa_y\zh\pa_y\hat{N}_0(\eh)-\zh^2\Delta(\hat{N}_0(\eh)).$$ The estimates \fref{boundnoone}, \fref{boundnotwo}, \fref{boundnothree}, \fref{boundnofour} still hold from direct check for $\hat{N}_0(\eh)$,
and we estimate with Lemma \ref{lemmainterpolation} and \fref{radiationbis}:
\bee
\int_{y\geq 1} |\hat{N}_0(\eh)H(\zh^2)|^2 & \lesssim & \int_{y\geq 1}\left[|\zh\Delta \zh|^2+\frac{\zh^4}{y^4}+|\pa_y\zh|^4\right]\\
& \lesssim &  \left[\left\|\frac{\eh}{y}\right\|_{L^{\infty}(y\geq 1)}^2+\left\|\pa_y\eh\right\|_{L^{\infty}(y\geq 1)}^2\right]\Sigma_{i=0}^2\int\frac{|\pa_y^i\eh|^2}{1+y^{2(3-i)}}\\
& \lesssim & b^2|\log b|^Cb^3|\log b|^C\lesssim b^5|\log b|^C.
\eee
\bee
\int_{y\geq 1} |\zh\pa_y\zh\pa_y\hat{N}_0(\eh)|^2 & \lesssim & |\log b|^C \left\|\frac{\eh}{y}\right\|_{L^{\infty}(y\geq 1)}^2\int_{y\geq 1}\left[\frac{1}{y^4}+b^2\right]\left[\frac{|\pa_y\eh|^2}{y^2}+\frac{|\eh|^2}{y^4}\right]\\
& \lesssim & b^2|\log b|^Cb^3|\log b|^C\lesssim b^5|\log b|^C,
\eee
$$
\int|\zh^2\Delta(\hat{N}_0(\eh))|^2\lesssim |\log b|^C \left\|\frac{\eh}{y}\right\|_{L^{\infty}(y\geq 1)}^2\int_{y\geq 1}\left[\frac{1}{y^6}+b^3\right]\frac{|\eh|^2}{y^2}\lesssim b^5|\log b|^C.$$
We thus conclude using \fref{htwoehat}: $$|(\eh_2,H\hat{N}(\eh))|\lesssim \left(b^2|\log b|^Cb^5\right)^{\frac12}\leq b^3|\log b|^2.$$
This concludes the proof of \fref{tobepreovehtwo}.\\
This concludes the proof of the Proposition \ref{bootstrap}.


\subsection{Proof of Theorem \ref{thmmain}}


We are now in position to conclude the proof of Theorem \ref{thmmain}. The proof relies on the reintegration of the modulation equations as in \cite{MR1}, \cite{MR4}, \cite{MRR}, we sketch the argument for the sake of completeness.\\

{\bf step 1} Finite time blow up.\\

Let $T\leq +\infty$ be the life time of the full map $v$ given by \fref{defvvu}, then the estimates of Proposition \ref{bootstrap} hold on $[0,T)$. From \fref{rougboundpope}, \fref{lawintegrationone}, $$-\frac{d}{dt}\sqrt\lambda=-\frac{1}{2\l\sqrt\l}\lsl\gtrsim \frac{b}{\l \sqrt{\lambda}}\gtrsim C(u_0)>0$$ and thus $\lambda$ touches zero at some finite time $T_0<+\infty$. Using \fref{formuleprojection}, \fref{formuleprojectionbis}, it is easily seen that the estimates of Proposition \ref{bootstrap} and the bootstrap bounds of Proposition \fref{lemmainterpolation} imply: $$\forall t\in [0,T_0), \ \ \|\Delta v(t)\|_{L^2}\lesssim C(t)<+\infty$$ and thus from the blow up criterion \fref{regularitiyt}: $$T_0=T<+\infty.$$ Observe then from \fref{lawintegrationone} that this implies 
\be
\label{boudnaryb}
\lambda(T)=b(T)=0.
\ee

{\bf step 2} Derivation of the sharp blow up speed.\\

We now slightly refine our control of $b$ through a logarithmic gain in the modulation equation \fref{parameterspresicely}. We commute \fref{eqepsilon}  with $H$ and take the inner product with $\chi_{B_\delta}\Lambda Q$ to derive:
\bea
\label{vnoooejeojier}
\nonumber &&\frac{d}{ds}\left\{(H\e,\chi_{B_\delta}\Lambda Q)\right\}-(H\e,\pa_s\chi_{B_\delta}\Lambda Q)+\lsl(\chi_{B_\delta}\Lambda Q,H\Lambda\e)+(H^2\e,\chi_{B_\delta}\Lambda Q)\\
& = & \left(H\left[-\Psit_b+L(\e)-N(\e)-Mod\right],\chi_{B_\delta}\Lambda Q\right).
\eea
We now estimate all terms in the above identity. First, for $\delta$ small enough, we estimate in brute force:
$$|(H\e,\pa_s\chi_{B_\delta}\Lambda Q)|+|\lsl(\chi_{B_\delta}\Lambda Q,H\Lambda\e)|+|(H[L(\e)-N(\e)],\chi_{B_\delta}\Lambda Q)|\lesssim \frac{b}{b^{C\delta}}\sqrt{\mathcal E_4}\lesssim \frac{b^2}{|\log b|^2}.$$ We then estimate the linear term: $$|(H^2\e,\chi_{B_\delta}\Lambda Q)|\lesssim \sqrt{\mathcal E_4}\sqrt{|\log b|}\lesssim \frac{b^2}{\sqrt{|\log b|}}.$$ The leading order $\Psit_b$ term is computed from \fref{descriptionerreur}, \fref{decomppsib}:
$$(-H\Psit_b,\chi_{B_{\delta}}\Psi_b)=-b^2(H\Sigma_b,\chi_{B_{\delta}}\Psi_b)+O\left(\frac{b^3}{b^{C\delta}}\right)=b^2c_b(\Lambda Q,\chi_{B_{\delta}}\Lambda Q)+O\left(\frac{b^2}{|\log b|^2}\right).$$ Finally, we compute the modulation term from \fref{defmode}:
\bee
(-HMod,\chi_{B_\delta}\Lambda Q) & = & \left(\lsl+b\right)(H\Lambda \qbt,\chi_{B_\delta}\Lambda Q)-(b_s+b^2)(H(\tt_1+2b\tt_2),\chi_{B_\delta}\Lambda Q)\\
& = & (b_s+b^2)(\Lambda Q,\chi_{B_\delta}\Lambda Q)+O\left(\frac{b}{b^{C\delta}}\frac{b^2}{|\log b|}\right).
\eee
We thus inject the collection of above estimates into \fref{vnoooejeojier} and derive the modulation equation:
$$(b_s+b^2)(\Lambda Q,\chi_{B_\delta}\Lambda Q)=\frac{d}{ds}\left\{(H\e,\chi_{B_\delta}\Lambda Q)\right\}-c_bb^2(\Lambda Q,\chi_{B_{\delta}}\Lambda Q)+O\left(\frac{b^2}{\sqrt{|\log b|}}\right)$$ which we rewrite using \fref{cb} and an integration by parts in time:
\bea
\label{intito}
&&\frac{d}{ds}\left\{b-\frac{(H\e,\chi_{B_\delta}\Lambda Q)}{(\Lambda Q,\chi_{B_\delta}\Lambda Q)}\right\}+b^2\left(1+\frac{2}{|\log b|}\right)\\
\nonumber & = & O\left(\frac{b^2}{|\log b|^{\frac32}}\right)+(H\e,\chi_{B_\delta}\Lambda Q)\frac{(\Lambda Q,\pa_s\chi_{B_\delta}\Lambda Q)}{(\Lambda Q,\chi_{B_\delta}\Lambda Q)^2}.
\eea
We now estimate: 
$$
\left|(H\e,\chi_{B_\delta}\Lambda Q)\frac{(\Lambda Q,\pa_s\chi_{B_\delta}\Lambda Q)}{(\Lambda Q,\chi_{B_\delta}\Lambda Q)^2}\right|\lesssim \frac{\sqrt{\mathcal E_4}}{b^{C\delta}}\frac{|b_s|}{b}\lesssim \frac{b^3}{b^{C\delta}},$$
$$\left|\frac{(H\e,\chi_{B_\delta}\Lambda Q)}{(\Lambda Q,\chi_{B_\delta}\Lambda Q)}\right|\lesssim \frac{\sqrt{\mathcal E_4}}{b^{C\delta}}\lesssim \frac{b^2}{b^{C\delta}}.
$$
We injecte these bounds into \fref{intito} and conclude that the quantity
\be
\label{lienbbtilde}
\tilde{b}=b-\frac{(H\e,\chi_{B_\delta}\Lambda Q)}{(\Lambda Q,\chi_{B_\delta}\Lambda Q)}=b+O\left(\frac{b^2}{|\log b|^2}\right)
\ee
satisfies the pointwise differential control:$$\left|\tilde{b}_s+\bt^2\left(1+\frac{2}{|\log \tilde{b}|}\right)\right|\lesssim \frac{\tilde{b}^2}{|\log \tilde{b}|^{\frac32}}.$$ Equivalently, $$\frac{\bt_s}{\bt^2\left(1+\frac{2}{|\log \bt|}\right)}+1=O\left(\frac{1}{|\log \bt|^{\frac32}}\right).$$ 
We now integrate this in time using $\lim_{s\to +\infty}\bt(s)=0$ from \fref{boudnaryb}, \fref{lienbbtilde} and get:
$$\bt(s)=\frac1s-\frac{2}{s\log s}+O\left(\frac{1}{s|\log s|^{\frac32}}\right)$$ and thus from \fref{lienbbtilde}: 
\be
\label{esitmateforb}
b(s)=\frac1s-\frac{2}{s\log s}+O\left(\frac{1}{s|\log s|^{\frac32}}\right).
\ee We now inject the modulation equation \fref{modulationequations} and conclude: $$-\lsl=\frac1s-\frac{2}{s\log s}+O\left(\frac{1}{s|\log s|^{\frac32}}\right).$$ 
We rewrite this as $$\left|\frac{d}{ds}\log\left(\frac{s\lambda(s)}{(\log s)^2}\right)\right|\lesssim \frac{1}{s|\log s|^{\frac32}}$$ and thus integrating in time yields the existence of $\kappa(u)>0$ such that: 
$$\frac{s\lambda(s)}{(\log s)^2}=\frac{1}{\kappa(u)}\left[1+O\left(\frac{1}{|\log s|^{\frac32}}\right)\right].
$$
Taking the log yields the bound $$|\log \lambda|=|\log s|\left[1+O\left(\frac{|\log \log s|}{\log s}\right)\right]$$ and thus $$\frac{1}{s}=\kappa(u)\frac{\lambda}{|\log \lambda|^2}\left(1+o(1)\right).$$ Injecting this into \fref{esitmateforb} yields:
\be
\label{cnoencoenoe}
-\lambda\lambda_t=-\lsl=\frac1s\left(1+o(1)\right)=\kappa(u)\frac{\lambda}{|\log \lambda|^2}\left(1+o(1)\right)
\ee and thus $$-|\log \lambda|^2\lambda_t=\kappa(u)(1+o(1)).$$ Integrating from $t$ to $T$ with $\lambda(T)=0$ yields 
$$\lambda(t)=\kappa(u)\frac{T-t}{|\log (T-t)|^2}\left[1+o(1)\right],$$ and \fref{universallawkgeq} is proved. This also implies using \fref{rougboundpope}: 
\be
\label{lineblambda}
b(t)=\kappa^2(u)\frac{T-t}{|\log (T-t)|^4}\left[1+o(1)\right].
\ee
In particular: 
\be
\label{realtionblb}
\frac{b}{\lambda}=\frac{k(u)}{|\log b|^2}(1+o(1)).
\ee
 
 {\bf step 3} $\dot{H}^2$ bound.\\
 
 We now turn to the proof of \fref{convustarb}, \fref{regularityustar}. Let $v(t,x)$ be the map associated to $u(t,r)$, explicitely: $$v(t,x)=\left|\begin{array}{lll}g(u(t,r)\cos\theta\\g(u(t,r))\sin\theta\\z(u(t,r)\end{array}\right.$$ and $\mathcal Q$ be given by \fref{defharmmap}. Let then 
 \be
 \label{defvtilde}
 \tilde{v}(t,x)=v(t,x)-\mathcal Q\left(\frac{x}{\lambda(t)}\right),
 \ee 
 and correspondingly 
 \be
 \label{dcnconeoneo}
 \tilde{u}(t,r)=u(t,r)-Q\left(\frac{r}{\lambda(t)}\right)=\left(\alphat+\e\right)_{\lambda(t)}.
 \ee
We claim the bound:
 \be
 \label{vounds}
 \forall t\in [0,T), \ \ \|\Delta \tilde{v}(t,x)\|_{L^2}\leq C(v_0).
 \ee
 Indeed, let the normal vector to the revolution surface $M$ at $v$ be given by $${\bf n}=\left|\begin{array}{lll}-z'(u)\cos\theta\\-z'(u)\sin\theta\\g'(u)\end{array}\right. $$ and compute the Laplace operator: $$\Delta v=\left[\Delta v-(\Delta v\cdot {\bf n}){\bf n}\right]+(\Delta v\cdot {\bf n}){\bf n}$$ with explicitely: 
 \be
 \label{formuleprojection}
 \Delta v-(\Delta v\cdot {\bf n}){\bf n}=\left(\Delta u-\frac{f(u)}{r^2}\right)\left|\begin{array}{lll} g'(u)\cos\theta\\g'(u)\sin\theta\\ z'(u)\end{array}\right.,
 \ee
 \be
 \label{formuleprojectionbis}\Delta v\cdot {\bf n}{\bf n}=\left|\begin{array}{lll}\left[(\pa_ru)^2g''(u)-\frac{g(u)(z'(u))^2}{r^2}\right]\cos\theta\\\left[(\pa_ru)^2g''(u)-\frac{g(u)(z'(u))^2}{r^2}\right]\sin\theta\\ (\pa_ru)^2z''(u)-\frac{z(u)(z'(u))^2}{r^2}\end{array}\right., 
 \ee we now claim: 
 \be
 \label{cneoeofheofreo}
 \int|\Delta \vt-(\Delta \vt\cdot {\bf n}){\bf n}|^2\lesssim C(v_0),
 \ee
 \be
 \label{cneoeofheofreobis}
 \int|\Delta \vt\cdot {\bf n}|^2\lesssim C(v_0),
 \ee
 which implies \fref{vounds}.\\
 {\it Proof of \fref{cneoeofheofreo}}: We inject the decomposition \fref{dcnconeoneo} and estimate:
 \bea
\label{cohneofihofe}
\nonumber && \int|\Delta \vt-(\Delta \vt\cdot {\bf n}){\bf n}|^2 = \int\left|\Delta u-\frac{f(u)}{r^2}\right|^2=\frac{1}{\l^2}\int \left|\Delta(\alphat+\e) -\frac{f(Q+\alphat+\e)-f(Q)}{y^2}\right|^2\\
 & \lesssim & \frac{1}{\l^2}\left[\int|H\e|^2+\int|H\alphat|^2+\int|f''(Q)|^2\frac{|\alphat|^4+|\e|^4}{y^4}+\int \frac{|\alphat|^6+|\e|^6}{y^4} \right].
 \eea
 We compute from the definition of $\alphat$ and \fref{realtionblb}:
 $$\int |H\alphat|^2\lesssim \int_{y\leq 2B_1}\left[by|\log y|+b^3\frac{y^3}{b|\log b|}\right|^2\lesssim b^2|\log b|^2\lesssim \l^2,$$ $$\int\frac{|\alphat|^4+|\alphat|^6}{y^4}\lesssim \int_{y\leq 2B_1}\frac{1}{y^4}\left[by|\log y|+b^3\frac{y^3}{b|\log b|}\right|^4\lesssim b^3|\log b|^C\lesssim \lambda^2.$$ Observe now from \fref{cnikheoeif},  \fref{universallawkgeq}, \fref{lineblambda} the bound:
$$\int|H\e(t)|^2= \mathcal E_2(t)\lesssim C(u_0)\lambda^2(t)\left[1+\int_0^T\frac{b^3(\tau)|\log b(\tau)|^2}{\lambda^4(\tau)}d\tau\right]\lesssim C(u_0)\lambda^2(t)$$ where we used from \fref{universallawkgeq}, \fref{lineblambda} again: $$\int_0^T\frac{b^3(\tau)|\log b(\tau)|^2}{\lambda^4(\tau)}d\tau\lesssim \int_0^T\frac{1}{(T-\tau)|\log (T-\tau)|^2}d\tau<+\infty.$$ We now estimate using \fref{estlinftydeux}:
$$\int_{y\leq 1} |f''(Q)|^2\frac{|\e|^4}{y^4}\lesssim b^6\lesssim \l^2,$$ and using $|f''(Q)|\lesssim \frac{y}{1+y^2}$, \fref{estlinftysecond} and the energy bound: $$ \int_{y\geq 1} |f''(Q)|^2\frac{|\e|^4}{y^4}\lesssim \left\|\frac{\e}{y^2}\right\|_{L^{\infty}(y\geq 1)}^2\int\frac{|\e|^2}{y^2}\lesssim b^3|\log b|^C\lesssim \l^2.$$ Injecting these bounds into \fref{cohneofihofe} yields \fref{cneoeofheofreo}.\\
{\it Proof of \fref{cneoeofheofreobis}}: We first claim:
\be
\label{fhohoehohoeone}
\int\left|(\pa_ru)^2g''(u)-(\pa_rQ_\lambda)^2g''(Q_\l)\right|^2\lesssim 1.
\ee
Indeed, we estimate:
\bee
\left|(\pa_ru)^2g''(u)-(\pa_rQ_\lambda)^2g''(Q_\l)\right| & \lesssim & |(\pa_ru)^2(g''(u)-g''(Q_\l))|+|g''(Q_\l)||\pa_ru|^2-|\pa_rQ_\l|^2|\\
& \lesssim & |\ut|\left[|\pa_rQ_\l|^2+|\pa_r\ut|^2\right]+|g''(Q_\l)|\left[|\pa_rQ_\l||\pa_r\ut|+|\pa_r\ut|^2\right]
\eee
and thus after rescaling using Lemma \ref{lemmainterpolation}:
\bee
&&\l^2\int\left|(\pa_ru)^2g''(u)-(\pa_rQ_\lambda)^2g''(Q_\l)\right|^2\lesssim \int\frac{|\alphat|^2+|\e|^2}{1+y^8}+\int|\pa_y(\alphat+\e)|^4|\alphat+\e|^2\\
&  + &  \int\frac{1}{1+y^6}|\pa_y(\alphat+\e)|^2+\|\pa_y(\alphat+\e)\|_{L^{\infty}}^2\int \frac{1}{1+y^2}|\pa_y(\alphat+\e)|^2\\
& \lesssim & b^2+\int(|\pa_y\alphat|^4+|\pa_y\e|^4)(|\alphat|^2+|\e|^2).
\eee
We estimate:
$$\int|\pa_y\alpha|^4|\alphat|^2\lesssim b^4|\log b|^C,\ \ \int|\pa_y\e|^4|\alphat|^2\lesssim \|\pa_y\e\|_{L^{\infty}}^2\|\alphat\|_{L^{\infty}}^2\lesssim b^4|\log b|^C,$$$$\int|\pa_y\alphat|^4|\e|^2\lesssim b^4|\log b|^C\int_{y\leq 2B_1}|\e|^2\lesssim b^4|\log b|^CB_1^8\int\frac{|\e|^2}{(1+y^8)|\log b|^C}\lesssim b^4|\log b|^C,$$  
and using a two dimensional Gagliardo Nirenberg inequality for the radial function $\pa_y\e$:
\be
\label{esttwodeirvtivlfour}
\int|\pa_y\e|^4|\e|^2\lesssim b^4+\|\e\|_{L^{\infty}}^2\|\pa_{yy} \e\|_{L^2}^2\|\pa_y\e\|_{L^2}^2\lesssim b^4+\mathcal E_2\lesssim \l^2,
\ee 
which concludes the proof of \fref{fhohoehohoeone}.\\
We now claim: 
\be
\label{novnoojoeroep}
\int\left|\frac{g(u)(z'(u))^2-g(Q_\l)(z'(Q_\l))^2}{r^2}\right|^2\lesssim 1.
\ee
Indeed,
\bee
\left|\frac{g(u)(z'(u))^2-g(Q_\l)(z'(Q_\l))^2}{r^2}\right| & \lesssim & \frac{|z'(Q_\l)|^2+|\ut|+g(Q_\l)}{r^2}|\ut|
\eee
and then from $z'(0)=z'(\pi)=0$:
$$\int \frac{|z'(Q)|^4+|g(Q)|^2}{y^4}|\alphat+\e|^2\lesssim \int\frac{|\alphat|^2+|\et|^2}{y^2(1+y^4)}\lesssim b^3|\log b|^C\lesssim \l^2,$$
\bee
\int \frac{|\ut|^2}{y^4}|\alphat+\e|^4 & \lesssim & b^3|\log b|^C+\int_{y\geq 1} \frac{|\e|^4}{y^4}\\
& \lesssim &  b^3|\log b|^C+\int\left[ |A\e|^4+|\pa_y\e|^4\right]\lesssim b^3 |\log b|^C+\l^2+\int|\nabla (A\e)|^2\int|A\e|^2\\
& \lesssim &\l^2+\int |A^*A\e|^2\lesssim \l^2
\eee
 where we used \fref{esttwodeirvtivlfour}, and \fref{novnoojoeroep} follows by rescaling. We estimate along the same lines:
 $$\int\left|(\pa_ru)^2z''(u)-(\pa_rQ_\l)^2z''(Q_\l)-\frac{z(u)(z'(u))^2-z(Q_\l)(z'(Q_\l))^2}{r^2}\right|^2\lesssim 1,$$ this is left to the reader.
 This concludes the proof of \fref{cneoeofheofreobis}\\
 
{\bf step 3} Quantization of the focused energy.\\

We now turn to the proof of \fref{convustarb}, \fref{regularityustar} and adapt the strategy in \cite{MR5}. The regularity of $v(t,x)$ outside the origin is a standard consequence of parabolic regularity and the fact that in corotational symmetry the nonlinearity is singular at the origin only. Hence there exists $v^*\in \dot{H}^1$ such that $$\forall R>0, \ \ \nabla u(t)\to \nabla v^*\ \ \mbox{in} \ \ L^2(|x|\geq R)\ \ \mbox{as}\ \ t\to T.$$ Moreover, $v$ is $\dot{H^1}$ bounded by conservation of energy, and thus recalling the decomposition \fref{defvtilde} and the uniform bound \fref{vounds}:
$$ \nabla \tilde{v}(t)\to \nabla v^*\ \ \mbox{in} \ \ L^2\ \ \mbox{and}\ \ \Delta v^*\in L^2$$ which concludes the proof of \fref{convustarb}, \fref{regularityustar}. This concludes the proof of Theorem \ref{thmmain}.


\begin{appendix}



\section{Regularity in corotational symmetry}


We detail in this appendix the regularity of $\dot{H}^{1}\cap \dot{H}^4$ maps with 1-corotational symmetry.

\begin{lemma}[Regularity in corotational symmetry]
\label{regularitycorot}
Let $v$ be a 1-corotational map with 
\be
\label{smallenness}
\|\nabla v-\nabla \mathcal Q\|_{\dot{H}^1}\ll 1 \ \ \mbox{and}\ \ \|v\|_{\dot{H}^i}<+\infty, \ \ 2\leq i\leq 4.
\ee Then $v$ admits a representation 
\be
\label{corotmap}
v(x)=\left|\begin{array}{lll} g(u(y))\cos\theta\\g(u(y))\sin\theta\\z(u(y))\end{array}\right.\ \ \mbox{where}\ \ u(y)=Q(y)+\e(y)
\ee 
satisfies the boundary conditions 
\be
\label{boudnaryconditions}
\e(0)=0, \ \ \lim_{y\to+\infty}\e(y)=0,
\ee the Sobolev bounds 
\be
\label{sobolefofiue}
\Sigma_{i=1}^4\int|\pa^{(i)}_y \e|^2+\int\frac{|\e|^2}{y^2}<+\infty
\ee
 and the regularity at the origin:
\bea
\label{chohohfoe}
\nonumber &&\int|H^2\e|^2+\int \frac{|H \e|^2}{y^4(1+|\log y|)^2}+ \int\frac{|AH \e|^2}{y^2}+\int|\pa_y(AH\e)|^2+\int \frac{|\pa^3_y \e|^2}{y^2(1+|\log y|)^2}\\
\nonumber &+&\int \frac{|\pa_y^2 \e|^2}{y^4(1+|\log y|)^2}+\int \frac{|\pa_y \e|^2}{y^2(1+|\log y|)^2}
+ \int \frac{|\e|^2}{y^4(1+|\log y|)^2}\\
& < & +\infty.
\eea
\end{lemma}

{\bf Proof of Lemma \ref{regularitycorot}} \\

{\bf step 1} $\dot{H}^1, \dot{H}^2$ bound.\\

From \fref{corotmap}, $$\int|\nabla v|^2=2\pi\int\left[|\pa_yu|^2+\frac{g^2(u)}{y^2}\right]$$ and thus the structure of $g$, the Sobolev embedding in radial symmetry 
\be
\label{sobiboej}
\|u\|^2_{L^{\infty}}\lesssim \|\pa_y u\|_{L^2}\|\frac{u}{y}\|_{L^2}
\ee and the smalness \fref{smallenness} ensure: $$\lim_{y \to 0}\e(y)=\lim_{y\to +\infty}\e(y)=0.$$ The energy bound 
\be
\label{energybound}
\int|\pa_y\e|^2+\frac{|\e|^2}{y^2}<+\infty
\ee
easily follows. Moreover, the energy density $$e(y)=|\pa_1v|^2+|\pa_2v|^2=2\left[|\pa_yu|^2+\frac{g^2(u)}{y^2}\right]$$ is bounded near the origin from the regularity $v\in \dot{H}^1\cap\dot{H}^4$ which implies 
\be
\label{boundu} 
|\pa_yu|+\frac{|u(y)|}{y}\lesssim 1\ \ \mbox{for} \ \ y\le1.
\ee
We now recall \fref{formuleprojection} which implies using the $\dot{H}^2$ boundedness of $v$: $$\int\left|\Delta u-\frac{f(u)}{y^2}\right|^2<+\infty.$$ Using \fref{energybound}, \fref{boundu} and the structure of $f$, we conclude: $$\int|Hu|^2<+\infty.$$  

{\bf step 2} $\dot{H}^4$ bound.\\

We now recall \fref{formuleprojection}, \fref{formuleprojectionbis} which we rewrite:
\be
\label{cioeoeo}
\Delta v=H\e\left|\begin{array}{lll}\cos\theta\\\sin\theta\\0\end{array}\right.+F(\e)
\ee
with:
\bee
F(\e)& = & \left(\Delta u-\frac{f(u)}{y^2}\right)\left|\begin{array}{lll} (g'(u)-1)\cos\theta\\(g'(u)-1)\sin\theta\\ z'(u)\end{array}\right.+\frac{f(Q+\e)-f(Q)-\e f'(Q)+(f'(Q)-1)\e}{y^2}\left|\begin{array}{lll}\cos\theta\\\sin\theta\\0\end{array}\right.\\
& + & \left|\begin{array}{lll}\left[(\pa_ru)^2g''(u)-\frac{g(u)(z'(u))^2}{r^2}\right]\cos\theta\\\left[(\pa_ru)^2g''(u)-\frac{g(u)(z'(u))^2}{r^2}\right]\sin\theta\\ (\pa_ru)^2z''(u)-\frac{z(u)(z'(u))^2}{r^2}\end{array}\right., 
 \eee
 We then compute $\nabla \Delta v$, $\Delta^2v$ and use the odd parity of $f,g$ and the cancellation at the origin \fref{boundu}  to conclude after a brute force computation: 
 \be
 \label{newboudnary}
 \int|\pa_yH\e|^2+\int\frac{|H\e|^2}{y^2}+\int|H^2\e|^2<+\infty.
 \ee The Sobolev bound \fref{sobolefofiue} away from the origin now easily follows: $$\Sigma_{i=1}^4\int_{y\geq 1}|\pa^{(i)}_y \e|^2<+\infty.$$
 
 {\bf step 3} Regularity at the origin.\\
 
 We claim that \fref{newboudnary} also implies the regularity \fref{chohohfoe} at the origin. Indeed, we estimate from \fref{sobiboej} and the two dimensional Hardy inequality with logarithmic loss: 
 \bea
 \label{pouetun}
\nonumber && \int_{y\leq 1}\left[\frac{|\pa_yH\e|^2}{y^2(1+|\log y|^2)}+\frac{|H\e|^2}{y^4(1+|\log y|^2)}\right] \\
\nonumber & \lesssim &  1+\int \frac{|\nabla \Delta v|^2}{y^2(1+|\log y|^2)}\lesssim 1+\int|\Delta^2v|^2+\int_{1\leq y\leq 2}|\nabla \Delta v|^2\\
&  \lesssim & 1
 \eea
 We similarily compute from \fref{formuleprojection} after an explicit computation: $$\int_{y\leq 1}\frac{(AH\e)^2}{y^2}\lesssim 1+\int\left|\frac{1}{r}\pa_{\theta}\pa_1\left(\Delta v-(\Delta v\cdot {\bf n}){\bf n}\right)\right|^2<+\infty$$

 We now observe from $$\pa_y(\log (\Lambda \phi))=\frac{Z}{y}$$ that for any function $h$: $$A^*h=\pa_yh+\frac{1+Z}{y}h=\frac{1}{y\Lambda \phi}\pa_y(y\Lambda \phi h),$$ and thus using the a priori bound \fref{boundu}: 
$$A\e(y)=\frac{1}{y\Lambda \phi(y)}\int_0^y \tau\Lambda\phi(\tau)H\e(\tau)d\tau.
$$
We then estimate from Cauchy Schwarz and Fubbini:
\bea
\label{pouetdeux}
\nonumber &&\int_{y\le 1}\frac{|A\e|^2}{y^5(1+|\log y|^2)}dy  \lesssim \int_{0\le y\leq 1}\int_{0\le\tau\le y}\frac{y^5}{y^9(1+|\log y|^2)}|H\e(\tau)|^2dyd\tau\\
\nonumber & \lesssim & \int_{0\le \tau\le 1}|H\e(\tau)|^2\left[\int_{\tau\leq y\leq 1}\frac{dy}{y^4(1+|\log y|^2)} \right]d\tau\lesssim  \int_{\tau\le 1}\frac{|H\e(\tau)|^2}{\tau^3(1+|\log \tau|^2)} d\tau\\
&\lesssim & 1.
\eea
We now rewrite near the origin: $$H\e=-\pa^2_y\e+\frac{1}{y}\left(-\pa_y \e+\frac{\e}{y}\right)+\frac{V-1}{y^2}\e=-\pa_y^2\e+\frac{A\e}{y}+\frac{(V-1)+(1-Z)}{y^2}\e$$ which implies using \fref{pouetun}, \fref{pouetdeux}, \fref{boundu}:
\bee
\int_{y\leq 1}\frac{|\pa^2_y\e|^2}{y^4(1+|\log y|^2)} & \lesssim & \int_{y\le 1}\frac{|H\e|^2}{y^4(1+|\log y|^2)}+\int_{y\le 1}\frac{|A\e|^2}{y^6(1+|\log y|^2)}+\int_{y\le 1}\frac{|\e|^2}{y^4(1+|\log y|^2)}\\
& \lesssim & 1,
\eee
\bee\nonumber \int_{y\leq 1}\frac{|\pa^3_y\e|^2}{y^2(1+|\log y|^2)}&\lesssim& \int_{y\le 1}\frac{|\pa_yH\e|^2}{y^2(1+|\log y|^2)}+ \int_{y\le 1}\frac{|H\e|^2}{y^4(1+|\log y|^2)}\\
& + & \int_{y\le 1}\frac{|\pa_y\e|^2}{y^2(1+|\log y|^2)}+\int_{y\le 1}\frac{|\e|^2}{y^4(1+|\log y|^2)}\\
& \lesssim & 1 
\eee
This concludes the proof of Lemma \ref{regularitycorot}.
 

\section{Coercivity bounds and interpolation estimates}


We recall in this section the coercivity bounds we use involving the operators $H,\tilde{H}$ and their iterate. Let us start with the coercivity of $\tilde{H}$:

\begin{lemma}[Coericivity of $\tilde{H}$]
\label{coerchtilde}
Let $\e_3$ with $$\int\frac{|\e_3|^2}{y^2}<+\infty,$$ then 
\be
\label{coercwthree}
\mathcal E_4=(\tilde{H}\e_3,\e_3)=\|A^*\e_3\|_{L^2}^2\geq c_0\left[\int|\pa_y\e_3|^2+\int\frac{|\e_3|^2}{y^2(1+|\log y|^2)}\right]
\ee
for some universal constant $c_0>0$.
\end{lemma}

{\bf Proof of Lemma \ref{coerchtilde}}. The proof is immediate in the case of $\Bbb S^2$ target thanks to the sign $$\tilde{V}=\frac{4}{1+y^2}\geq 0$$ which yields \fref{coercwthree} from standard two dimensional weigthed Hardy inequality -see \cite{MRR} for more details-. For a general $g$ however, this sign property is lost and the proof relies on a standard compactness argument and the explicit knowledge of the kernel of $A^*$.\\
We first claim the subcoercivity property: 
\bea
\label{subcoerc}
\nonumber (\tilde{H}\e_3,\e_3)&\geq& c_0\left(\int|\pa_y\e_3|^2+\int_{y\leq 1}\frac{|\e_3|^2}{y^2}+\int_{y\geq 1}\frac{\e_3^2}{y^2(1+y^2)}+\int\frac{|\e_3|^2}{y^2(1+|\log y|^2)}\right)\\
& - &  \frac{1}{c_0}\int\frac{|\e_3|^2}{1+y^4}
\eea
for some universal constant $c_0>0$. Indeed, we estimate near the origin using the expansion \fref{comportementvtilde}: $$\int_{y\leq1 }\frac{\tilde{V}}{y^2}|\e_3|^2\geq c_0\int_{y\leq 1}\frac{|\e_3|^2}{y^2}-\frac{1}{c_0}\int_{y\leq 1}|\e_3|^2.$$ For $y\geq 1$, the logarithmic Hardy bound  $$\int\frac{|\e_3|^2}{y^2(1+|\log y|^2)}\lesssim \int|\nabla\e_3|^2+\int_{1\leq y\leq 2}|\e_3|^2$$ and the degeracy from \fref{comportementvtilde} $|\tilde{V}(y)|\lesssim \frac{1}{y^2}$ yield \fref{subcoerc}.\\
We now prove \fref{coercwthree} and argue by contradiction. Let a sequence $\e_3^{(n)}$ with 
\be
\label{infordec}
\int|\pa_y\e^{(n)}_3|^2+\int_{y\leq 1}\frac{|\e^{(n)}_3|^2}{y^2}+\int_{y\geq 1}\frac{|\e^{(n)}_3|^2}{y^2(1+|\log y|^2)}=1, \ \ (\tilde{H}\e_3^{(n)},\e_3^{(n)})\leq \frac1n,
\ee then the sequence $\e_3^{(n)}$ is bounded in $H^1_{loc}$ and thus weakly converges to $\e_3^*$ up to a subsequence. Moreoever, \fref{subcoerc}, $$\int\frac{|\e_3^{(n)}|^2}{1+y^4}\gtrsim c_0>0$$ and thus by lower semi continuity of norms and the compactness of the Sobolev embedding $H^1\hookrightarrow L^2_{loc}$: 
\be
\label{orifini}
\int|\pa_y\e^*_3|^2+\int_{y\leq 1}\frac{|\e^*_3|^2}{y^2}+\int_{y\geq 1}\frac{|\e^*_3|^2}{y^2(1+|\log y|^2)}=1, \ \ \int\frac{|\e_3^{*}|^2}{1+y^4}\gtrsim c_0>0
\ee 
and $$(\tilde{H}\e_3^*,\e_3^*)=\|A^*\e_3^*\|_{L^2}^2\leq 0.$$ Hence $$A^*\e_3=0\ \ \mbox{ie} \ \ \e_3=\frac{c}{y\Lambda Q}$$ which contradicts the boundary condition at the origin imposed by \fref{orifini}. This concludes the proof of \fref{coercwthree} and Lemma \ref{coerchtilde}.\\

We now claim the following coercivity properties for $H,H^2$ which rely on a similar compactness argument and the explicit knowledge of the kernel of $H,H^2$. The proof is given in \cite{RR}, \cite{MRR} in the case $g(u)=\sin u$, but the same argument applies under our assumptions on $g$, the key being the behaviour at the origin and infinity \fref{comportementz}, \fref{comportementv}, \fref{comportementvtilde}. The proof is therefore left to the reader:

\begin{lemma}[Coercivity of $H,H^2$]
\label{lemmacoer}
Let $M>0$ and $\e$ radially symmetric satisfying \fref{boudnaryconditions}, \fref{sobolefofiue}, \fref{chohohfoe} and the orthogonality conditions $$(\e,\Phi_M)=(\e,H\Phi_M)=0.$$ Then there exists $C(M)>0$ such that the following bounds hold:
$$(H\e,\e)=\int|A\e|^2\geq C(M)\int\left[(\pa_y\e)^2+\frac{\e^2}{y^2}\right],$$
$$\int(H\e)^2\geq C(M)\int \left[\frac{(\pa_y\e)^2}{y^2(1+|\log y|^2)}+\frac{\e^2}{y^4(1+|\log y|^2)}\right],$$
\bee
\int(H^2\e)^2& \geq & C(M)\int\left[\frac{|H \e|^2}{y^4(1+|\log y|)^2}+ \frac{|\pa_y H \e|^2}{y^2(1+|\log y|)^2}+\int \frac{|\pa^3_y \e|^2}{y^2(1+|\log y|)^2}\right.\\
& +&\left .\frac{|\pa_y^2 \e|^2}{y^4(1+|\log y|)^2}+\int \frac{|\pa_y \e|^2}{y^2(1+y^4)(1+|\log y|)^2}
+ \int \frac{|\e|^2}{y^4(1+y^4)(1+|\log y|)^2}\right]
\eee
\end{lemma}

We now recall the interpolation bounds on $\e$ needed all along the proof of Proposition \ref{bootstrap} which are a direct consequence of two dimensional weighted Hardy estimates and the coercivity bounds of Lemma \ref{lemmacoer}. These bounds which hold in the setting of the bootstrap bounds \fref{init1h}, \fref{init2h}, \fref{init3h}, \fref{init3hbis} were explicitely derived in \cite{MRR} to which we refer for a proof.

\begin{lemma}[Interpolation estimates]
\label{lemmainterpolation}
There holds -with constants a priori depending on M-:
\bea
\label{estun}
 & & \int \frac{|\e|^2}{y^4(1+y^4)(1+|\log y|^2)}+\int \frac{|\pa^i_y\e|^2}{y^2(1+y^{6-2i})(1+|\log y|^2)}\\
\nonumber& & \lesssim \mathcal E_4,\ \ 1\leq i\leq 3,
\eea
\be
\label{lossyboundwperp}
\int_{y\geq 1} \frac{1+|\log y|^C}{y^2(1+|\log y|^2)(1+y^{6-2i})}|\pa_y^i\e|^2 \lesssim b^4|\log b|^{C_1(C)}, \ \ 0\leq i\leq 3,
\ee
\be
\label{inteproloatedbound}
\int_{y\geq 1} \frac{1+|\log y|^C}{y^2(1+|\log y|^2)(1+y^{4-2i})}|\pa_y^i\e|^2 \lesssim b^3|\log b|^{C_1(C)}, \ \ 0\leq i\leq 3,
\ee
\be
\label{estperplinfty}
\|\e\|_{L^{\infty}}\lesssim \delta(\alpha^*),
\ee
\be
\label{estinterm}
\|A\e\|^2_{L^{\infty}}\lesssim b^2|\log b|^2,
\ee
\be
\label{estoriiginagian}
\left\|\frac{A \e}{y^2(1+|\log y|)}\right\|_{L^{\infty}(y\leq 1)}^2+\left\|\frac{\Delta A \e}{1+|\log y|}\right\|^2_{L^{\infty}(y\leq 1)}+\left\|\frac{H \e}{y(1+|\log y|)}\right\|_{L^{\infty}(y\leq 1)}^2\lesssim b^4,
\ee
\be
\label{choeouefeouie}
\left\|\frac{|H \e|}{y(1+|\log y|)}\right\|_{L^{\infty}(y\leq 1)}^2\lesssim b^4,
\ee
\be
\label{estlinftydeux}
\|\frac{\e}{y}\|^2_{L^{\infty}(y\leq 1)}+\left\|\frac{\pa_y\e}{\sqrt{1+|\log y|}}\right\|^2_{L^{\infty}(y\leq 1)}\lesssim b^4,
\ee
\be
\label{estlinftydeuxfa}
\|\frac{\e}{y}\|^2_{L^{\infty}(y\geq 1)}+\|\pa_y\e\|^2_{L^{\infty}(y\geq 1)}\lesssim b^2|\log b|^C,
\ee
\be
\label{estlinftysecond}
\|\frac{\e}{1+y^2}\|^2_{L^{\infty}}+\|\frac{\pa_y\e}{1+y}\|^2_{L^{\infty}}+\|\pa_{yy}\e\|^2_{L^{\infty}(y\geq 1)}\lesssim b^3|\log b|^C,
\ee
\be
\label{estoeofj}
\|\pa_{yyy}\e\|_{L^{\infty}(y\geq 1)}^2\lesssim b^4|\log b|^C.
\ee
\end{lemma}

\end{appendix}


\end{document}